\documentclass[12pt]{article}
\usepackage[psamsfonts]{amssymb}

\title{\bf Aggregation functions for decision making\footnote{Chapter 17 of {\sl Concepts and methods of decision making} (ISTE/John Wiley).}}

\author{Jean-Luc Marichal \\
\small{Mathematics Research Unit, University of Luxembourg} \\
\small{162A, avenue de la Fa\"{\i}encerie, L-1511 Luxembourg} \\
\small{jean-luc.marichal[at]uni.lu}}

\date{January 27, 2009}

\begin{document}
\maketitle

\newtheorem{definition}{Definition}[section]
\newtheorem{lemma}{Lemma}[section]
\newtheorem{proposition}{Proposition}[section]
\newtheorem{theorem}{Theorem}[section]
\newtheorem{corollary}{Corollary}[section]

\newcommand{\N}{\mathbb{N}}                     
\newcommand{\R}{\mathbb{R}}                     
\newcommand{\Vspace}{\vspace{2ex}}              
\newcommand{\fn}{{\cal F}_N}

\section{Introduction}

Aggregation functions are generally defined and used to combine several numerical values into a single one, so that the final
result of the aggregation takes into account all the individual values in a given manner. Such functions are widely used in many well-known
disciplines such as statistics, economics, finance, and computer science. For general background, see Grabisch et al.~\cite{GraMarMesPap09}.

For instance, suppose that several individuals form quantifiable judgements either about a measure of an object (weight, length, area, height,
volume, importance or other attributes) or about a ratio of two such measures (how much heavier, longer, larger, taller, more important,
preferable, more meritorious etc. one object is than another). In order to reach a consensus on these judgements, classical aggregation
functions have been proposed: arithmetic mean, geometric mean, median and many others.

In multicriteria decision making, values to be aggregated are typically {\sl preference}\/ or {\sl satisfaction}\/ degrees. A preference degree
reveals to what extent an alternative $a$ is preferred to an alternative $b$, and thus is a relative appraisal. By contrast, a satisfaction degree
expresses to what extent a given alternative is satisfactory with respect to a given criterion. It is an absolute appraisal.

We assume that the values to be aggregated belong to numerical scales, which can be of ordinal or cardinal type. On an ordinal scale, numbers
have no meaning other than defining an order relation on the scale; distances or differences between values cannot be interpreted. On a
cardinal scale, distances between values are not quite arbitrary. There are actually several kinds of cardinal scales. On an interval scale,
where the position of the zero is a matter of convention, values are defined up to a positive linear transformation i.e. $\phi(x)=rx+s$, with
$r>0$ and $s\in\R$ (e.g.\ temperatures expressed on the Celsius scale). On a ratio scale, where a true zero exists, values are defined up to a
similarity transformation i.e. $\phi(x)=rx$, with $r>0$ (e.g.\ lengths expressed in inches). We will come back on these measurement aspects in
Section~\ref{sec:pslte}.

Once values are defined we can aggregate them and obtain a new value. This can be done in many different ways according to what is expected
from the aggregation function, the nature of the values to be aggregated, and which scale types have been used. Thus, for a given
problem, any aggregation function should not be used. In other terms, the use of a given aggregation function should always be justified.

To help the practitioner choose an appropriate aggregation function in a given problem, it is useful and even convenient to adopt an axiomatic
approach. Such an approach consists in classifying and choosing aggregation functions according to the properties they fulfill. Thus, a catalog
of ``desirable'' properties is proposed and, whenever possible, a description of the family of aggregation functions satisfying a given set of
properties is provided. This is the very principle of axiomatization.

Proposing an interesting axiomatic characterization of an aggregation function (or a family of aggregation functions) is not an easy task.
Mostly, aggregation functions can be characterized by different sets of conditions. Nevertheless the various possible characterizations are not
equally important. Some of them involve purely technical conditions with no clear interpretation and the result becomes useless. Others
involve conditions that contain the result explicitly and the characterization becomes trivial. On the contrary, there are characterizations
involving only natural conditions which are easily interpretable. In fact, this is the only case where the result should be seen as a significant
contribution. It improves our understanding of the function considered and provides strong arguments to justify (or reject) its use in a given
context.

The main aim of this chapter is to present, on an axiomatic basis, the most used families of aggregation functions in decision making. We shall
confine ourselves to aggregation functions that assign a numerical value to every profile of $n$ values, which represent objects or
alternatives. We will not deal with utility functions which, in a more general way, make it possible to rank alternatives without assigning
precise values to them. For instance, procedures such as `leximin' or `discrimin' are ranking procedures, rather than aggregation functions.

The organization of this chapter is as follows. In Section 2 we yield the list of the main properties that we shall use. This list is divided
into three classes: (1) elementary properties (continuity, symmetry, etc.); (2) properties related to the scale types used to represent the
data; and (3) certain algebraic properties such as associativity. In Section 3 we present the concept of mean and its various definitions. Perhaps
the most common definition of means is that of quasi-arithmetic means with a very natural axiomatization due to Kolmogoroff and Nagumo. In
Section 4 we present associative functions, which are at the root of the theory of semi-groups. These functions permitted to develop the concept
of fuzzy connectives such as t-norms, t-conorms, and uninorms. In Section 5 we present an important branch of the aggregation function theory,
namely Choquet and Sugeno non-additive integrals. These integrals enable us to generalize the classical aggregation modes, such as the weighted
arithmetic mean and the median, to functions that take into account the possible interactions among the considered attributes. Finally, in
Sections 6 and 7 we present particular functions designed for aggregating interval scales, ratio scales, and ordinal scales.

We close this introduction by setting the notation that we will use in this chapter.

In a general manner, we shall denote an aggregation function with $n$ variables by $A:E^n\to \R$ where $E$ is a real interval, bounded or not.
$E^{\circ}$ will denote the interior of $E$. We shall sometimes consider sequences of functions $(A^{(n)}:E^n\to\R)_{n\ge 1}$, the superscript
$^{(n)}$ being used only to specify the number of arguments of the function $A^{(n)}$.

We shall use $N$ to denote the index set $\{1,\ldots,n\}$ and $2^N$ to denote the set of its subsets. $\Pi_N$ will be used to denote the set of
permutations on $N$. Finally, for any $S\subseteq N$, the characteristic vector of $S$ in $\{0,1\}^n$ will be denoted $\mathbf{1}_S$.

There is also a rather standard notation for certain aggregation functions. Here are the most common ones:

\begin{itemize}
\item The {\sl arithmetic mean}\/ is defined as $${\rm AM}(x)=\frac 1n\sum_{i=1}^n x_i.$$%
\item For any weight vector $\omega=(\omega_1,\ldots,\omega_n)\in [0,1]^n$ such that $\sum_i \omega_i=1$, the {\sl weighted arithmetic mean}\/
and the {\sl ordered weighted averaging function}\/ are defined as
\begin{eqnarray*}
{\rm WAM}_{\omega}(x)&=&\sum_{i=1}^n\omega_ix_i,\\
{\rm OWA}_{\omega}(x)&=&\sum_{i=1}^n\omega_ix_{(i)},
\end{eqnarray*}
respectively, where $(\cdot)$ represents a permutation on $N$ such that $x_{(1)}\le\cdots\le x_{(n)}$.%
\item For any $k\in N$, the {\sl projection}\/ and the {\sl order statistic} associated with the $k$th argument are defined as
\begin{eqnarray*}
{\rm P}_k(x) &=& x_k, \\
{\rm OS}_k(x) &=& x_{(k)},
\end{eqnarray*}
respectively.

\item For any $S\subseteq N$, $S\neq\varnothing$, the {\sl partial minimum}\/ and {\sl partial maximum}\/ functions associated with $S$ are defined as
\begin{eqnarray*}
{\rm min}_S(x) &=& \min_{i\in S}x_i,\\
{\rm max}_S(x) &=& \max_{i\in S}x_i,\\
\end{eqnarray*}
respectively.
\end{itemize}

In this chapter the min and max operations will often be denoted $\wedge$ and $\vee$, respectively.

\section{Aggregation properties}

As mentioned in the introduction, in order to choose a reasonable or satisfactory aggregation mode, it is useful to adopt an axiomatic approach
and impose that the aggregation functions fulfill some selected properties. Such properties can be dictated by the nature of the values to be
aggregated. For example, in some multicriteria evaluation methods, the aim is to assess a global absolute score to an alternative given a set of
partial scores with respect to different criteria. Clearly, it would be unnatural to give as a global score a value which is lower than the
lowest partial score, or greater than the highest score, so that only internal aggregation functions (means) are allowed. Another example
concerns the aggregation of opinions in voting procedures. If, as usual, the voters are anonymous, the aggregation function must be symmetric.

In this section we present some properties that could be desirable or not depending upon the considered problem. Of course, all these
properties are not required with the same strength, and do not pertain to the same purpose. Some of them are imperative conditions whose
violation leads to obviously counterintuitive aggregation modes. Others are technical conditions that simply facilitate the representation or the
calculation of the aggregation function. There are also facultative conditions that naturally apply in special circumstances but are not to be
universally accepted.

\subsection{Elementary mathematical properties}

\begin{definition}
$A:E^n\to\R$ is {\sl symmetric}\/ if, for any $\pi\in\Pi_N$, we have
$$A(x_1,\ldots,x_n)=A(x_{\pi(1)},\ldots,x_{\pi(n)}) \qquad (x\in E^n).$$
\end{definition}

The symmetry property essentially implies that the indexing (ordering) of the arguments does not matter. This is required when combining
criteria of equal importance or the opinions of anonymous experts.

\begin{definition}
$A:E^n\to\R$ is {\sl continuous}\/ if it is continuous in the usual sense.
\end{definition}

One of the advantages of a continuous aggregation function is that it does not present any chaotic reaction to a small change of the arguments.

\begin{definition}
$A:E^n\to\R$ is
\begin{itemize}
\item {\sl nondecreasing}\/ if, for any $x,x'\in E^n$, we have $$x\le x' \quad\Rightarrow\quad A(x)\le A(x'),$$ %
\item {\sl strictly increasing}\/ if it is nondecreasing and if, for any $x,x'\in E^n$, we have
$$x\le x'
~\mbox{et}~ x\neq x' \quad\Rightarrow\quad A(x)< A(x'),$$%
\item{\sl unanimously increasing}\/ if it is nondecreasing and if, for any $x,x'\in E^n$, we have
$$x<x' \quad\Rightarrow\quad A(x)< A(x').$$
\end{itemize}
\end{definition}

An increasing aggregation function presents a non-negative response to any increase of the arguments. In other terms, increasing a partial value
cannot decrease the result. This function is strictly increasing if, moreover, it presents a positive reaction to any increase of at least one
argument. Finally, a unanimously increasing function is increasing and presents a positive response whenever all the arguments strictly
increase. For instance we observe that on $[0,1]^n$, the maximum function $A(x)=\max x_i$ is unanimously increasing whereas the bounded sum
$A(x)=\min(\sum_{i=1}^nx_i,1)$ is not.

\begin{definition}
$A:E^n\to\R$ is {\sl idempotent}\/ if $A(x,\ldots,x)=x$ for all $x\in E$.
\end{definition}

\begin{definition}
$A:[a,b]^n\to\R$ is {\sl weakly idempotent}\/ if $A(a,\ldots,a)=a$ and $A(b,\ldots,b)=b$.
\end{definition}

In a variety of applications, it is desirable that the aggregation functions satisfy the idempotency property, i.e.\ if all $x_i$ are identical,
$A(x_1,\ldots,x_n)$ restitutes the common value.

\begin{definition}\label{def:cdi}
$A:E^n\to\R$ is
\begin{itemize}
\item {\sl conjunctive}\/ if $A(x)\le\min x_i$ for all $x\in E^n$,%
\item {\sl disjunctive}\/ if $\max x_i\le A(x)$ for all $x\in E^n$,%
\item {\sl internal}\/ if $\min x_i\le A(x)\le \max x_i$ for all $x\in E^n$.
\end{itemize}
\end{definition}

Conjunctive functions combine values as if they were related by a logical {\bf AND} operator. That is, the result of aggregation can be high
only if all the values are high. t-norms are suitable functions for doing conjunctive aggregation (see Section~\ref{sec:ttu}). At the opposite,
disjunctive functions combine values as an {\bf OR} operator, so that the result of aggregation is high if at least one value is high. The best
known disjunctive functions are t-conorms.

Between these two extreme situations are the internal functions, located between the minimum and the maximum of the arguments. In this kind of
functions, a bad (resp.\ good) score on one criterion can be compensated by a good (resp.\ bad) one on another criterion, so that the result of
aggregation will be medium. By definition, means are internal functions (see Section~\ref{sec:moy}).

\subsection{Stability properties related to scale types}
\label{sec:pslte}

Depending on the kind of scale which is used, allowed operations on values are restricted. For example, aggregation on ordinal scales should be
limited to operations involving comparisons only, such as medians and order statistics.

A {\sl scale of measurement}\/ is a mapping which assigns real numbers to objects being measured. The {\sl type}\/ of a scale, as defined by
Stevens~\cite{Ste51,Ste59}, is defined by a class of {\sl admissible transformations}, transformations that lead from one acceptable scale to
another. For instance, we call a scale a {\sl ratio scale}\/ if the class of admissible transformations consists of the similarities
$\phi(x)=rx$, with $r>0$. In this case, the scale value is determined up to choice of a unit. Mass is an example of a ratio scale. The
transformation from kilograms into pounds, for example, involves the admissible transformation $\phi(x)=2.2x$. Length (inches, centimeters) and
time intervals (years, seconds) are two other examples of ratio scales. We call a scale an {\sl interval scale}\/ if the class of admissible
transformations consists of the positive linear transformations $\phi(x)=rx+s$, with $r>0$ and $s\in\R$. The scale value is then determined up
to choices of unit and zero point. Temperature (except where there is an absolute zero) defines an interval scale. Thus, transformation from
Centigrade into Fahrenheit involves the admissible transformation $\phi(x)=9x/5+32$. We call a scale an {\sl ordinal scale}\/ if the class of
admissible transformations consists of the strictly increasing bijections $\phi(x)$. Here the scale value is determined only up to order. For
example, the scale of air quality being used in a number of cities is an ordinal scale. It assigns a number 1 to unhealthy air, 2 to
unsatisfactory air, 3 to acceptable air, 4 to good air, and 5 to excellent air. We could just as well use the numbers 1, 7, 8, 15, 23, or the
numbers 1.2, 6.5, 8.7, 205.6, 750, or any numbers that preserve the order. Definitions of other scale types can be found in the book by Roberts
\cite{Rob79} on measurement theory, see also Roberts~\cite{Rob90,Rob94}. The reader will find further details on measurement in
Chapter~18 of the present volume.

A statement using scales of measurement is said to be {\sl meaningful}\/ if the truth or falsity of the statement is invariant when every scale
is replaced by another acceptable version of it \cite[p.~59]{Rob79}. For example, a ranking method is meaningful if the ranking of alternatives
induced by the aggregation does not depend on scale transformation.

In 1959, Luce~\cite{Luc59} observed that the general form of a functional relationship between variables is greatly restricted if we know the
scale type of the variables. These restrictions are discovered by formulating a functional equation from knowledge of the admissible
transformations. Luce's method is based on the principle of theory construction, which states that an admissible transformation of the
independent variables may lead to an admissible transformation of the dependent variable. For example, suppose that
$f(a)=A(f_1(a),\ldots,f_n(a))$, where $f$ and $f_1,\ldots,f_n$ are all ratio scales, with the units chosen independently. Then, by the principle
of theory construction, we obtain the functional equation
$$\displaylines{A(r_1x_1,\ldots,r_nx_r)=R(r_1,\ldots,r_n)A(x_1,\ldots,x_n), \cr r_i>0, ~~~ R(r_1,\ldots,r_n)>0.}$$

Acz\'el et al.~\cite{AczRobRos86} showed that the solutions of this equation are given by
$$A(x)=a\prod_{i=1}^ng_i(x_i), ~~~\mbox{with $a>0$, $g_i>0$,}$$ and $$g_i(x_iy_i)=g_i(x_i)g_i(y_i).$$

In this section we present some functional equations related to certain scale types. The interested reader can find more details in Acz\'el et
al.~\cite{AczRob89,AczRobRos86} and a good survey in Roberts~\cite{Rob94}.

\begin{definition}\label{de:signi}
$A:\R^n\to\R$ is
\begin{itemize}
\item {\sl meaningful for the same input-output ratio scales}\/
if, for any $r>0$, we have $$A(rx_1,\ldots,rx_n)=rA(x_1,\ldots,x_n)\qquad (x\in\R^n),$$%
\item {\sl meaningful for the same input ratio scales}\/
if, for any $r>0$, there exists $R_r>0$ such that $$A(rx_1,\ldots,rx_n)=R_rA(x_1,\ldots,x_n)\qquad (x\in\R^n),$$%
\item {\sl meaningful for the same input-output interval scales}\/
if, for any $r>0$ and $s\in\R$, we have $$A(rx_1+s,\ldots,rx_n+s)=rA(x_1,\ldots,x_n)+s\qquad (x\in\R^n),$$%
\item {\sl meaningful for the same input interval scales}\/ if, for any $r>0$ and $s\in\R$, there exist $R_{r,s}>0$ and $S_{r,s}\in\R$ such that
$$A(rx_1+s,\ldots,rx_n+s)=R_{r,s}A(x_1,\ldots,x_n)+S_{r,s}\qquad (x\in\R^n),$$%
\item {\sl meaningful for the same input-output ordinal scales}\/ if, for any strictly increasing bijection $\phi:\R\to\R$,
we have $$A(\phi(x_1),\ldots,\phi(x_n))=\phi(A(x_1,\ldots,x_n))\qquad (x\in\R^n),$$%
\item {\sl meaningful for the same input ordinal scales}\/ if, for any strictly increasing bijection $\phi:\R\to\R$, there exists a strictly
increasing function $\psi_{\phi}:\R\to\R$ such that
$$A(\phi(x_1),\ldots,\phi(x_n))=\psi_{\phi}(A(x_1,\ldots,x_n))\qquad (x\in\R^n).$$%
\end{itemize}
\end{definition}

\subsection{Algebraic properties}

The following properties concern the aggregation procedures that can be decomposed into partial aggregations, that is, procedures for which it
is possible to partition the set of attributes into disjoint subgroups, build the partial aggregation for each subgroup, and then combine these
partial results to get the global value. This condition may take several forms. Maybe one of the strongest is associativity. Other weaker
formulations will also be presented: decomposability and bisymmetry.

We first present associativity for two variable functions.

\begin{definition}\label{eq:asso2}
$A:E^2\to E$ is {\sl associative}\/ if, for any $x\in E^3$, we have
$$A(A(x_1,x_2),x_3)=A(x_1,A(x_2,x_3)).$$
\end{definition}

A large number of papers deal with the associativity functional equation. For a list of references see Acz\'el~\cite[\S6.2]{Acz66}.

This property can be extended to sequences of functions as follows.

\begin{definition}
The sequence $(A^{(n)}:\R^n\to\R)_{n\ge 1}$ is {\sl associative}\/ if $A^{(1)}(x)=x$ for all $x\in E$ and
$$
A^{(n)}(x_1,\ldots,x_k,x_{k+1},\ldots,x_n)=A^{(n)}(A^{(k)}(x_1,\ldots,x_k),A^{(n-k)}(x_{k+1},\ldots,x_n))
$$
for all $x \in E^n$ and all $k,n \in \N$ such that $1\le k\le n$.
\end{definition}

Implicit in the assumption of associativity is a consistent way of going unambiguously from the aggregation of $n$ elements to $n+1$ elements,
i.e., if M is associative
$$A^{(n+1)}(x_1,\ldots,x_{n+1})=A^{(2)}(A^{(n)}(x_1,\ldots,x_n),x_{n+1}),$$ for all $n\in\N\setminus\{0\}$.

Let us turn to the decomposability property. For this purpose, we introduce the following notation: For any $k\in \N\setminus\{0\}$ and any
$x\in\R$, we set $k\cdot x=x,\ldots,x$ ($k$ times). For example,
$$A(3\cdot x,2\cdot y)=A(x,x,x,y,y).$$
\begin{definition}\label{def:dec}
The sequence $(A^{(n)}:\R^n\to\R)_{n\ge 1}$ is {\sl decomposable}\/ if $A^{(1)}(x)=x$ for all $x\in E$ and
$$
A^{(n)}(x_1,\ldots,x_k,x_{k+1},\ldots,x_n)=A^{(n)}(k\cdot A^{(k)}(x_1,\ldots,x_k),(n-k)\cdot A^{(n-k)}(x_{k+1},\ldots,x_n))
$$
for all $x \in E^n$ and all $k,n \in \N$ such that $1\le k\le n$.
\end{definition}

Here the definition is the same as for associativity, except that the partial aggregations are duplicated by the number of aggregated
values.

Introduced first in Bemporad~\cite[p.~87]{Bem26} in a characterization of the arithmetic mean, decomposability has been used by
Kolmogoroff~\cite{Kol30} and Nagumo~\cite{Nag30} to characterize the quasi-arithmetic means. More recently, Marichal et Roubens~\cite{MarRou93}
proposed calling this property ``decomposability" in order to avoid confusion with classical associativity.

The {\sl bisymmetry}\/ property, which extends associativity and symmetry simultaneously, is defined for $n$-variables functions as follows.
\begin{definition}
$A:E^n\to E$ is {\sl bisymmetric}\/ if
\begin{eqnarray*}
&   & A(A(x_{11},\ldots,x_{1n}),\ldots,A(x_{n1},\ldots,x_{nn}))\\ 
& = & A(A(x_{11},\ldots,x_{n1}),\ldots,A(x_{1n},\ldots,x_{nn}))
\end{eqnarray*}
for all square matrix $(x_{ij})\in E^{n \times n}$.
\end{definition}

For 2-variable functions, this property has been investigated from the algebraic point of view by using it mostly in structures without the
property of associativity; see Acz\'el~\cite[\S6.4]{Acz66} and Acz\'el and Dhombres~\cite[Chapter~17]{AczDho89}.

For a sequence of functions, this property becomes as described in the following definition.
\begin{definition}
The sequence $(A^{(n)}:\R^n\to\R)_{n\ge 1}$ is {\sl bisymmetric} if $A^{(1)}(x)=x$ for all $x\in E$ and
\begin{eqnarray*}
&   & A^{(p)}(A^{(n)}(x_{11},\ldots,x_{1n}),\ldots,A^{(n)}(x_{n1},\ldots,x_{pn}))\\ 
& = & A^{(n)}(A^{(p)}(x_{11},\ldots,x_{p1}),\ldots,A^{(p)}(x_{1n},\ldots,x_{pn}))
\end{eqnarray*}
for all $n,p\in\N\setminus\{0\}$ and all matrix $(x_{ij})\in E^{p \times n}$.
\end{definition}

\section{Means}
\label{sec:moy}

It would be very unnatural to propose a chapter on aggregation functions without dealing somehow with {\sl means}. Already discovered and
studied by the ancient Greeks (see for instance Antoine~\cite[Chapter~3]{Ant98}) the concept of mean has given rise today to a very wide field
of investigation with a huge variety of applications. Actually, a tremendous amount of literature on the properties of several means (such as
the arithmetic mean, the geometric mean, etc.) has already been produced, especially since the 19th century, and is still developing today. For
a good overview, see the expository paper by Frosini~\cite{Fro87} and the remarkable monograph by Bullen et al.~\cite{BulMitVas88}.

The first modern definition of mean was probably due to Cauchy~\cite{Cau21} who considered in 1821 a mean as an internal (Definition~\ref{def:cdi}) function.

The concept of mean as a {\sl numerical equalizer}\/ is usually ascribed to Chisini~\cite[p.~108]{Chi29}, who provided the following
definition:
\begin{quote}
Let $y = g(x_1,\ldots,x_n)$ be a function of $n$ independent variables $x_1,\ldots,x_n$ representing homogeneous quantities. A mean of
$x_1,\ldots,x_n$ with respect to the function $g$ is a number $M$ such that, if each of $x_1,\ldots,x_n$ is replaced with $M$, the function
value is unchanged, that is,
$$
g(M,\ldots,M) = g(x_1,\ldots,x_n).
$$
\end{quote}

When $g$ is considered as the sum, the product, the sum of squares, the sum of inverses, or the sum of exponentials, the solution of Chisini's
equation corresponds to the arithmetic mean, the geometric mean, the quadratic mean, the harmonic mean, and the exponential mean, respectively.

Unfortunately, as noted by de Finetti~\cite[p.~378]{deF31} in 1931, Chisini's definition is so general that it does not even imply that the
``mean'' (provided there exists a real and unique solution to Chisini's equation) fulfills Cauchy's internality property.

The following quote from Ricci~\cite[p.~39]{Ric15} could be considered as another possible criticism to Chisini's view.
\begin{quote}
... when all values become equal, the mean equals any of them too. The inverse proposition is not true. If a function of several variables takes
their common value when all variables coincide, this is not sufficient evidence for calling it a mean. For example, the function
$$ g(x_1,x_2,\ldots,x_n) = x_n + (x_n - x_1) + (x_n - x_2) + \cdots + (x_n - x_{n-1})
$$
equals $x_n$ when $x_1 = \cdots = x_n$, but it is even greater than $x_n$ as long as $x_n$ is greater than every other variable.
\end{quote}

In 1930, Kolmogoroff~\cite{Kol30} and Nagumo~\cite{Nag30} considered that the mean should be more than just a Cauchy mean or a numerical
equalizer. They defined a {\sl mean value}\/ to be a decomposable sequence of continuous, symmetric, strictly increasing (in each variable), and
idempotent real functions
$$
M^{(1)}(x_1) = x_1, M^{(2)}(x_1,x_2), \ldots, M^{(n)}(x_1,\ldots,x_n), \ldots.
$$
They proved, independently of each other, that these conditions are necessary and sufficient for the quasi-arithmeticity of the mean, that is,
for the existence of a continuous strictly monotonic function $f$ such that $M^{(n)}$ may be written in the form
\begin{equation}\label{eq:qam}
M^{(n)}(x_1,\ldots,x_n) = f^{-1} \Bigl[\, \frac 1n \sum_{i=1}^n f(x_i) \Bigr]
\end{equation}
for all $n\in\N\setminus\{0\}$.

The {\sl quasi-arithmetic means}\/ (\ref{eq:qam}) comprise most of the algebraic means of common use; see Table~\ref{tab:means}. However, some
means, such as the median, do not belong to this family.

\begin{table}[ht]
\begin{center}
$$
\begin{array}{|ccccc|}
\hline
f(x)    && M^{(n)}(x_1,\ldots,x_n)          && \mbox{name}      \\
\hline
    &&                      &&          \\
x   && \frac 1n \sum\limits_{i=1}^n x_i         && \mbox{arithmetic}    \\
    &&                      &&          \\
x^2 && \Big(\frac 1n \sum\limits_{i=1}^n x_i^2\Big)^{1/2}        && \mbox{quadratic} \\
    &&                      &&          \\
\log x  && \Big(\prod\limits_{i=1}^n x_i\Big)^{1/n}      && \mbox{geometric} \\
    &&                      &&          \\
x^{-1}  && \frac 1{\frac 1n \sum\limits_{i=1}^n \frac 1{x_i}}   && \mbox{harmonic}\\
    &&                      &&          \\
x^{\alpha} \; (\alpha \in \R\setminus\{0\})
    && \Big(\frac 1n \sum\limits_{i=1}^n x_i^{\alpha} \Big)^{1/\alpha}
                            && \mbox{root-mean-power} \\
    &&                      &&          \\
e^{\alpha \, x} \; (\alpha \in \R\setminus\{0\})
    && \frac 1{\alpha} \ln \Big(\frac 1n \sum\limits_{i=1}^n e^{\alpha \, x_i}\Big)
                            && \mbox{exponential}   \\
    &&                      &&          \\
\hline
\end{array}
$$
\caption{Examples of quasi-arithmetic means} \label{tab:means}
\end{center}
\end{table}

The above properties defining a mean value seem to be natural enough. For instance, one can readily see that, for increasing means, the
idempotency property is equivalent to Cauchy's internality, and both are accepted by all statisticians as requisites for means.

The decomposability property of means is rather natural. Under idempotency, this condition becomes equivalent to
$$
\displaylines{M^{(k)}(x_1,\ldots,x_k) = M^{(k)}(x'_1,\ldots,x'_k) \cr \Downarrow \cr M^{(n)}(x_1,\ldots,x_k,x_{k+1},\ldots,x_n) =
M^{(n)}(x'_1,\ldots,x'_k,x_{k+1},\ldots,x_n)}
$$
which states that the mean does not change when altering some values without modifying their partial mean.

The purpose of this section is not to present a state of the art of all the known results in this vast realm of means. Instead, we just skim the
surface of the subject by pointing out characterization results for the most-often used and best-known families of means.

The medians and, more generally, the order statistics (which are particular means designed to aggregate ordinal values) will be briefly
presented in Section~\ref{sec:agrord}.

\subsection{Quasi-arithmetic means}

As we just mentioned, quasi-arithmetic means were introduced from a very natural axiomatization. In this section, we investigate those means
both as $n$-variable functions and as sequences of functions. Results on this class of means can also be found in Bullen et
al.~\cite[Chapitre~4]{BulMitVas88}.

It was proved by Acz\'el~\cite{Acz48} (see also \cite[\S6.4]{Acz66} and \cite[Chapter~17]{AczDho89}) that the quasi-arithmetic means are the
only symmetric, continuous, strictly increasing, idempotent, real functions $M : E^n \to E$ which satisfy the bisymmetry condition. The
statement of this result is formulated as follows.

\begin{theorem}\label{thm:sy/co/sin/id/b}
$M : E^n \to E$ is a symmetric, continuous, strictly increasing, idempotent, and bisymmetric function if and only if there exists a continuous
and strictly monotonic function $f : E \to \R$ such that
\begin{equation}\label{eq:sy/co/sin/id/b}
M(x) = f^{-1} \Big[\, \frac 1 n \sum_{i = 1}^n f(x_i)\Big] \qquad (x \in E^n).
\end{equation}
\end{theorem}

The {\sl quasi-arithmetic means}\/ (\ref{eq:sy/co/sin/id/b}) are internal aggregation functions and cover a wide spectrum of means including
arithmetic, quadratic, geometric, harmonic; see Table~\ref{tab:means}.

The function $f$ occurring in equation (\ref{eq:sy/co/sin/id/b}) is called a {\sl generator}\/ of $M$. It was also proved that $f$ is determined up to a
linear transformation: With $f(x)$, every function $$g(x) = r f(x) + s \qquad (r,s \in \R, r \neq 0)$$ belongs to the same $M$, but no other
function.

In addition to Acz\'el's result, we also recall Kolmogoroff-Nagumo's result.

\begin{theorem}
The sequence $(M^{(n)}:E^n\to E)_{n\ge 1}$ is a decomposable sequence of symmetric, continuous, strictly increasing, and idempotent functions if
and only if there is a continuous and strictly monotonic function $f : E \to \R$ such that
$$
M^{(n)}(x) = f^{-1} \Big[\, \frac 1 n \sum_{i = 1}^n f(x_i)\Big] \qquad (x \in E^n).
$$
\end{theorem}

Nagumo~\cite{Nag30} investigated some subfamilies of the class of quasi-arithmetic means. He proved the following result (see also
\cite[\S4]{AczAls87} and \cite[Chapitre~15]{AczDho89}).

\begin{proposition}
Assume $E = \,]0,\infty[$ or a subinterval.

\noindent $(i)$ $M : E^n \to E$ is a quasi-arithmetic mean that is meaningful for the same input-output ratio scales if and only if
\begin{itemize}
\item either $M$ is the geometric mean:
$$
M(x) = \Bigl(\, \prod_{i=1}^n x_i\Bigr)^{\frac 1n} \qquad (x \in E^n),
$$
\item or $M$ is a root-mean-power: there exists $\alpha \in \R\setminus\{0\}$ such that
\begin{equation}\label{eq:root/power/mean}
M(x) = \Bigl(\, \frac 1n \sum_{i=1}^n x_i^{\alpha} \Bigr)^{\frac 1{\alpha}} \qquad (x \in E^n).
\end{equation}
\end{itemize}
$(ii)$ $M : E^n \to E$ is a quasi-arithmetic mean that is meaningful for the same input-output interval scales if and only if $M$ is the
arithmetic mean.
\end{proposition}

Let us denote the root-mean-power (\ref{eq:root/power/mean}) generated by $\alpha \in \R\setminus\{0\}$ by $M_{(\alpha)}$. It is well known
\cite[\S16]{BecBel65} that, if $\alpha_1 < \alpha_2$ then $M_{(\alpha_1)}(x) \le M_{(\alpha_2)}(x)$ for all $x \in \, ]0,+\infty[^n$ (equality
if and only if all $x_i$ are equal).

The family of root-mean-powers was studied by Dujmovi\'c~\cite{Duj74,Duj75} and then by Dyckhoff et Pedrycz~\cite{DycPed84}. It encompasses most
of traditionally known means: the arithmetic mean $M_{(1)}$, the harmonic mean $M_{(-1)}$, the quadratic mean $M_{(2)}$, and three limiting
cases: the geometric mean $M_{(0)}$, the minimum $M_{(-\infty)}$ and the maximum $M_{(+\infty)}$ (see, e.g., Abramowitz and
Stegun~\cite{AbrSte64}).

Let us return to Theorem~\ref{thm:sy/co/sin/id/b}. Note that Acz\'el~\cite{Acz48} also investigated the case where symmetry and idempotency are dropped
(see also \cite[\S6.4]{Acz66} and \cite[Chapitre~17]{AczDho89}). He obtained the following result.

\begin{theorem}\label{thm:co/sin/b}
$(i)$ $M : E^n \to E$ is a continuous, strictly increasing, idempotent, and bisymmetric function if and only if there exists a continuous and
strictly monotonic function $f : E \to \R$ and real numbers $\omega_1,\ldots,\omega_n > 0$ fulfilling $\sum_i \omega_i = 1$ such that
\begin{equation}\label{eq:co/sin/id/b}
M(x) = f^{-1} \Bigl[\, \sum_{i = 1}^n \omega_i \, f(x_i)\Bigr] \qquad (x \in E^n).
\end{equation}
$(ii)$ $M : E^n \to E$ is a continuous, strictly increasing, and bisymmetric function if and only if there exists a continuous and strictly
monotonic function $f : E \to \R$ and real numbers $p_1,\ldots,p_n > 0$ and $q \in \R$ such that
\begin{equation}\label{eq:co/sin/b}
M(x) = f^{-1} \Bigl[\, \sum_{i = 1}^n p_i \, f(x_i) + q \Bigr] \qquad (x \in E^n).
\end{equation}
\end{theorem}

The {\sl quasi-linear means}\/ (\ref{eq:co/sin/id/b}) and the {\sl quasi-linear functions}\/ (\ref{eq:co/sin/b}) are weighted aggregation
functions. The question of uniqueness with respect to $f$ is dealt with in detail in Acz\'el~\cite[\S6.4]{Acz66}. Table~\ref{tab:co/sin/id/b}
provides some special cases of quasi-linear means.

\begin{table}[htb]
\begin{center}
$$
\begin{array}{|ccccc|}
\hline
f(x)    && M(x)                 && \mbox{name of weighted mean}     \\
\hline
    &&                      &&          \\
x   && \sum\limits_{i=1}^n \omega_i \, x_i      && \mbox{arithmetic}    \\
    &&                      &&          \\
x^2 && \Bigl(\sum\limits_{i=1}^n \omega_i\, x_i^2\Bigr)^{1/2}  && \mbox{quadratic} \\
    &&                      &&          \\
\log x  && \prod\limits_{i=1}^n x_i^{\omega_i}      && \mbox{geometric} \\
    &&                      &&          \\
x^{\alpha} \; (\alpha \in \R\setminus\{0\})   && \Bigl(\sum\limits_{i=1}^n \omega_i\, x_i^{\alpha} \Bigr)^{1/\alpha}   && \mbox{root-mean-power}    \\
    &&                      &&          \\
\hline
\end{array}
$$
\caption{Examples of quasi-linear means} \label{tab:co/sin/id/b}
\end{center}
\end{table}

\subsection{Lagrangian and Cauchy means}

Let us consider the intermediate point $M$ in the classical mean value formula
\begin{equation}\label{eq:mvt}
F(y)-F(x)=F'(M)(y-x) \qquad (x,y\in E),
\end{equation}
as a function of the variables $x,y$, with the convention $M(x,x)=x$, where $F:E\to \R$ is a given continuously differentiable and strictly
convex or strictly concave function. Reformulating this definition in terms of integrals instead of derivatives, we can rewrite (\ref{eq:mvt})
as
\begin{equation}\label{eq:lm}
M(x,y)=\cases{\displaystyle{f^{-1}\left( {1\over y-x}\int_x^y f(\xi) d\xi\right)}, & if $x\neq y$,\cr x, & if $x=y$,}
\end{equation}
for $x,y\in I$, where $f:E\to\R$ is a continuous strictly monotonic function. This function $M(x,y)$ is called the {\sl Lagrangian mean}\/
associated with $f$. See for example Berrone and Moro~\cite{BerMor98} and Bullen et al.~\cite[p.~343]{BulMitVas88}. The uniqueness of the generator
is the same as for quasi-arithmetic means, that is, defined up to a linear transformation; see Berrone and Moro~\cite[Corollary 7]{BerMor98} and
Matkowski~\cite[Theorem 1]{Mat99}.

Many classical means are Lagrangian. The arithmetic mean and the geometric means correspond to taking $f(x)=x$ and $f(x)=1/x^2$, respectively,
in (\ref{eq:lm}). The harmonic mean, however, is not Lagrangian.

In general, some of the most common means are both quasi-arithmetic and Lagrangian. However, there are quasi-arithmetic means, such as the harmonic one,
which are not Lagrangian. Conversely, the logarithmic mean
$$M(x,y)=\cases{\displaystyle{\frac{x-y}{\log x-\log y}}, & for $x,y>0$, $x\neq y$,\cr x, &
for $x=y>0$,}$$ is an example of a Lagrangian mean ($f(x)=1/x$) that is not quasi-arithmetic.

Let us now consider the Cauchy mean value theorem, which asserts that, for any functions $F$ and $g$, continuous on an interval $[x,y]$ and
differentiable on $]x,y[$, there exists $M\in\left]a,b\right[$ such that
$$\frac{F(y)-F(x)}{g(y)-g(x)}=\frac{F'(M)}{g'(M)}$$ If the functions $g$ and $f:=F'/g'$ are strictly
monotonic on $]x,y[$, the mean value $M(x,y)$ is unique and can be written as
$$M(x,y)=\cases{\displaystyle{f^{-1}\left( {1\over g(y)-g(x)}\int_x^y f(\xi) dg(\xi)\right)}, & if $x\neq y$,\cr x, & if
$x=y$,}$$ for $x,y\in E$. It is then said to be the {\sl Cauchy mean associated with the pair $(f,g)$}; see Berrone and Moro~\cite{BerMor00}.
Such a mean is continuous, idempotent, symmetric, and strictly increasing.

When $g=f$ (resp.\ $g$ is the identity function), we retrieve the quasi-arithmetic (resp.\ the Lagrangian) mean generated by $f$. The {\sl
anti-Lagrangian mean}\/ \cite{BerMor00} is obtained when $f$ is the identity function. For example, the harmonic mean is an anti-Lagrangian mean
generated by the function $g=1/x^2$. The generators of the same anti-Lagrangian mean are defined up to the same non-zero affine transformation.

\section{Associative aggregation functions}
\label{sec:conjdisj}

Before dealing with associative functions and their axiomatizations, we will need to introduce some useful concepts. A {\sl semigroup}\/ $(E,
A)$ is a set $E$ with an associative operation $A:E^2\to E$ defined on it. As usual, we assume that $E$ is a real interval, bounded or not.

An element $e \in E$ is\par $a)$ an {\sl identity}\/ for $A$ if $A(e,x) = A(x,e) = x$ for all $x \in E$, \par b) a {\sl zero}\/ (or {\sl
annihilator}) for $A$ if $A(e,x) = A(x,e) = e$ for all $x \in E$,\par c) an {\sl idempotent}\/ for $A$ if $A(e,e) = e$.

\Vspace

\noindent For any semigroup $(E, A)$, it is clear that there is at most one identity and at most one zero for $A$ in $E$, and both are
idempotents.

We also need to introduce the concept of {\sl ordinal sum}, well known in the theory of semigroups (see, e.g., Climescu~\cite{Cli46} and
Ling~\cite{Lin65}).

\begin{definition}\label{de:ord/sum}
Let $K$ be a totally ordered set and $\{(E_k, A_k) \, | \, k \in K\}$ be a collection of disjoint semigroups indexed by $K$. Then the {\sl
ordinal sum}\/ of $\{(E_k, A_k) \, | \, k \in K\}$ is the set-theoretic union $\cup_{k \in K} E_k$ under the following binary operation:
$$
A(x,y) = \cases{A_k(x,y),       & si $\exists \, k \in K$ such that $x,y \in E_k$ \cr
        \min(x,y),  & si $\exists \, k_1, k_2 \in K, k_1 \neq k_2$ such that $x \in E_{k_1}$ and $y \in E_{k_2}$.\cr}
$$
\end{definition}
The ordinal sum is a semigroup under the above defined operation.

\subsection{Strictly increasing functions}

Acz\'el~\cite{Acz48b} investigated the general continuous, strictly increasing, real solutions on $E^2$ of the associativity functional equation
(\ref{eq:asso2}). He proved the following result (see also \cite[\S6.2]{Acz66}).

\begin{theorem}\label{thm:co/sin/a2}
Let $E$ be a real interval, bounded or not, which is open on one side. $A:E^2\to E$ is continuous, strictly increasing, and associative if and
only if there exists a continuous and strictly monotonic function $f : E \to \R$ such that
\begin{equation}\label{eq:co/sin/a2}
A(x,y) = f^{-1}[f(x) + f(y)] \qquad ((x,y) \in E^2).
\end{equation}
\end{theorem}

It was also proved that the function $f$ occurring in (\ref{eq:co/sin/a2}) is determined up to a multiplicative constant, that is, with $f(x)$
all functions $g(x) = r \, f(x)$ ($r \in \R\setminus\{0\}$) belongs to the same $A$, and only these.

Moreover, the function $f$ is such that, if $e \in E$ then
\begin{equation}\label{eq:f/M/id}
A(e,e) = e ~~\Leftrightarrow ~~ f(e) = 0.
\end{equation}

By (\ref{eq:f/M/id}) and because of strict monotonicity of $f$, there is at most one idempotent for $A$ (which is, actually, the identity) and
hence $A$ cannot be idempotent. Therefore, there is no continuous, strictly increasing, idempotent, and associative function. However, we can
notice that every continuous, strictly increasing, and associative function is necessarily symmetric. The sum ($f(x) = x$) and the product
($f(x) = \log x$) are well-known examples of continuous, strictly increasing, and associative functions.

According to Ling~\cite{Lin65}, any semigroup $(E, M)$ satisfying the hypotheses of Theorem~\ref{thm:co/sin/a2} is called {\sl Acz\'elian}.

Recall that each associative sequence $(A^{(n)}:E^n\to E)_{n\ge 1}$ of functions is uniquely determined by its 2-variable function. We therefore have the following result.

\begin{corollary}\label{cor:co/sin/a}
Let $E$ be a real interval, bounded or not, which is open on one side. $(A^{(n)}:E^n\to E)_{n\ge 1}$ is an associative sequence of continuous
and strictly increasing functions if and only if there exists a continuous and strictly monotonic function $f : E \to \R$ such that, for all $n
\in \N\setminus\{0\}$,
$$
A^{(n)}(x) = f^{-1}\Bigl[\, \sum_{i=1}^n f(x_i)\Bigr] \qquad (x \in E^n).
$$
\end{corollary}

\subsection{Archimedean semigroups}

Some authors attempted to generalize Theorem~\ref{thm:co/sin/a2} by relaxing the strict increasing monotonicity into nondecreasing monotonicity.
However it seems that the class of continuous, nondecreasing, and associative functions has not yet been described. However, under some additional
conditions, results have been obtained.

First, we state a representation theorem attributed very often to Ling~\cite{Lin65}. In fact, her main theorem can be deduced from previously
known results on topological semigroups, see Faucett~\cite{Fau55} and Mostert and Shields~\cite{MosShi57}. Nevertheless, the advantage of Ling's
approach is twofold: it treats two different cases in a unified manner and establishes elementary proofs.

\begin{theorem}\label{thm:co/in/a2/b}
Let $E = [a,b]$. $A:E^2\to E$ is continuous, nondecreasing, associative, and
\begin{eqnarray}
&& A(b,x) = x \qquad (x \in E) \label{eq:mbx/x}\\
&& A(x,x) < x \qquad (x \in E^\circ) \label{eq:mxx/ix}
\end{eqnarray}
if and only if there exists a continuous strictly decreasing function $f : E \to [0,+\infty]$, with $f(b) = 0$, such that
\begin{equation}\label{eq:co/in/a2/b}
A(x,y) = f^{-1}[ \, \min(f(x) + f(y),f(a)) \, ] \qquad (x,y \in E).
\end{equation}
\end{theorem}

The requirement that $E$ be closed is not really a restriction. If $E$ is any real interval, finite or infinite with right endpoint $b$ ($b$
can be $+\infty$), then we can replace condition (\ref{eq:mbx/x}) with
$$
\lim_{t \to b^-} A(t,t) = b, \quad \lim_{t \to b^-} A(t,x) = x \qquad (x \in E).
$$

Any function $f$ solving equation (\ref{eq:co/in/a2/b}) is called an {\sl additive generator} (or simply {\sl generator}) of $A$. Moreover, we
can easily see that any function $A$ of the form (\ref{eq:co/in/a2/b}) is symmetric and conjunctive.

Condition (\ref{eq:mbx/x}) expresses that $b$ is a {\sl left identity} for $A$. It turns out, from (\ref{eq:co/in/a2/b}), that $b$ acts as an
identity and $a$ as a zero. Condition (\ref{eq:mxx/ix}) simply expresses that there are no idempotents for $A$ in $]a,b[$. Indeed, by
nondecreasing monotonicity and (\ref{eq:mbx/x}), we always have $A(x,x) \le A(b,x) = x$ for all $x \in [a,b]$.

Depending on whether $f(a)$ is finite or infinite (recall that $f(a) \in [0,+\infty]$), $A$ takes a well-defined form (see Fodor and
Roubens~\cite[\S1.3]{FodRou94} and Schweizer and Sklar~\cite{SchSkl83}):
\begin{itemize}
\item $f(a) < +\infty$ if and only if $A$ has {\em zero divisors} (i.e.\ $\exists \, x,y \in\, ]a,b[$ such that $A(x,y) = a$). In this case,
there exists a continuous strictly increasing function $g : [a,b] \to [0,1]$, with $g(a) = 0$ and $g(b) = 1$, such that
\begin{equation}\label{eq:luka-like}
A(x,y) = g^{-1}[\max(g(x) + g(y) - 1,0)] \qquad (x,y \in [a,b]).
\end{equation}
To see this, it suffices to set $g(x) := 1 - f(x)/f(a)$.

For associative sequences $(A^{(n)}:[a,b]^n\to [a,b])_{n\ge 1}$, (\ref{eq:luka-like}) becomes
$$
A^{(n)}(x) = g^{-1}\Bigl[ \max \Bigl(\sum_{i=1}^n g(x_i) - n + 1, 0\Bigr)\Bigr] \qquad (x \in [a,b]^n, ~ n \in \N\setminus\{0\}).
$$
\item $\lim_{t \to a^+} f(x) = +\infty$ if and only if $A$ is strictly increasing on $]a,b[$. In this case, there exists a continuous strictly
increasing function $g : [a,b] \to [0,1]$, with $g(a) = 0$ and $g(b) = 1$, such that
\begin{equation}\label{eq:prod-like}
A(x,y) = g^{-1}[ g(x) \, g(y)] \qquad (x,y \in [a,b]).
\end{equation}
To see this, it suffices to set $g(x) := \exp(-f(x))$.

For associative sequences $(A^{(n)}:[a,b]^n\to [a,b])_{n\ge 1}$, (\ref{eq:prod-like}) becomes
$$
A^{(n)}(x) = g^{-1}\Bigl[\, \prod_{i=1}^n g(x_i) \Bigr] \qquad (x \in [a,b]^n, ~ n \in \N\setminus\{0\}).
$$
\end{itemize}
Of course, Theorem~\ref{thm:co/in/a2/b} can also be written in a dual form as follows.

\begin{theorem}\label{thm:co/in/a2/a}
Let $E = [a,b]$. $A:E^2\to E$ is continuous, nondecreasing, associative, and
\begin{eqnarray*}
&& A(a,x) = x \qquad (x \in E)\\ 
&& A(x,x) > x \qquad (x \in E^\circ)
\end{eqnarray*}
if and only if there exists a continuous strictly increasing function $f : E \to [0,+\infty]$, with $f(a) = 0$, such that
\begin{equation}\label{eq:co/in/a2/a}
A(x,y) = f^{-1}[ \, \min(f(x) + f(y),f(b)) \, ] \qquad (x,y \in E).
\end{equation}
\end{theorem}

Again, $E$ can be any real interval, even infinite. Functions $A$ of the form (\ref{eq:co/in/a2/a}) are symmetric and disjunctive. There
are no interior idempotents. The left endpoint $a$ acts as an identity and the right endpoint $b$ acts as a zero.

Once more, two mutually exclusive cases can be examined:
\begin{itemize}
\item $f(b) < +\infty$ if and only if $A$ has zero divisors (i.e.\ $\exists \, x,y \in\, ]a,b[$ such that $A(x,y) = b$). In this case, there
exists a continuous strictly increasing function $g : [a,b] \to [0,1]$, with $g(a) = 0$ and $g(b) = 1$, such that
\begin{equation}\label{eq:coluka-like}
A(x,y) = g^{-1}[\min(g(x) + g(y),1)] \qquad (x,y \in [a,b]).
\end{equation}
To see this, it suffices to set $g(x) := f(x)/f(b)$.

For associative sequences $(A^{(n)}:[a,b]^n\to [a,b])_{n\ge 1}$, (\ref{eq:coluka-like}) becomes
$$
A^{(n)}(x) = g^{-1}\Bigl[ \min \Bigl(\sum_{i=1}^n g(x_i), 1\Bigr)\Bigr] \qquad (x \in [a,b]^n, ~ n \in \N\setminus\{0\}).
$$
\item $\lim_{t \to b^-} f(x) = +\infty$ if and only if $A$ is strictly increasing on $]a,b[$. In this case, there exists a continuous strictly
increasing function $g : [a,b] \to [0,1]$, with $g(a) = 0$ and $g(b) = 1$, such that
\begin{equation}\label{eq:coprod-like}
A(x,y) = g^{-1}[ 1 - (1 - g(x)) \, (1 - g(y))] \qquad (x,y \in [a,b]),
\end{equation}
To see this, it suffices to set $g(x) := 1 - \exp(-f(x))$.

For associative sequences $(A^{(n)}:[a,b]^n\to [a,b])_{n\ge 1}$, (\ref{eq:coprod-like}) becomes
$$
A^{(n)}(x) = g^{-1}\Bigl[ 1 - \prod_{i=1}^n (1 - g(x_i)) \Bigr] \qquad (x \in [a,b]^n, ~ n \in \N\setminus\{0\}).
$$
\end{itemize}

Any semigroup fulfilling the assumptions of Theorem~\ref{thm:co/in/a2/b} or \ref{thm:co/in/a2/a} is called {\sl Archime\-dean}; see
Ling~\cite{Lin65}. In other words, any semigroup $(E, A)$ is said to be {\sl Archimedean} if $A$ is continuous, nondecreasing, associative, one
endpoint of $E$ is an identity for $A$, and there are no idempotents for $A$ in $E^\circ$. We can make a distinction between conjunctive and
disjunctive Archimedean semigroups depending on whether the identity is the right or left endpoint of $E$, respectively. An
Archimedean semigroup is called {\sl properly Archimedean} or {\sl Acz\'elian} if every additive generator $f$ is unbounded; otherwise it is
{\sl improperly Archimedean}.

Ling~\cite[\S6]{Lin65} proved that every Archimedean semigroup is obtainable as a limit of Acz\'elian semigroups.

\subsection{A class of nondecreasing and associative functions}

We now give a description of the class of functions $A:[a,b]^2\to [0,1]$ that are continuous, nondecreasing, weakly idempotent, and associative.
For all $\theta \in [a,b]$, we define ${\cal A}_{a,b,\theta}$ as the set of continuous, nondecreasing, weakly idempotent, and associative
functions $A:[a,b]^2\to [0,1]$ such that $A(a,b) = A(b,a) = \theta$. The extreme cases ${\cal A}_{a,b,a}$ and ${\cal A}_{a,b,b}$ will play an
important role in the sequel. The results proved by the author can be found in \cite{Mar00b}.

\begin{theorem}\label{thm:co/in/wid/a}
$A:[a,b]^2\to [0,1]$ is continuous, nondecreasing, weakly idempotent, and associative if and only if there exist $\alpha, \beta \in [a,b]$ and
two functions $A_{a,\alpha \wedge \beta, \alpha \wedge \beta} \in {\cal A}_{a,\alpha \wedge \beta, \alpha \wedge \beta}$ and $A_{\alpha \vee
\beta, b, \alpha \vee \beta} \in {\cal A}_{\alpha \vee \beta, b, \alpha \vee \beta}$ such that, for all $x,y \in [a,b]$,
$$
A(x,y) = \cases{A_{a,\alpha \wedge \beta, \alpha \wedge \beta}(x,y), & if $x,y \in [a,\alpha \wedge \beta]$ \cr
        A_{\alpha \vee \beta, b, \alpha \vee \beta}(x,y), & if $x,y \in [\alpha \vee \beta,b]$ \cr
        (\alpha \wedge x) \vee (\beta \wedge y) \vee (x \wedge y), & otherwise.\cr}
$$
\end{theorem}

Now, let us turn to the description of ${\cal A}_{a,b,a}$. Mostert and Shields~\cite[p.\ 130, Theorem B]{MosShi57} proved the following

\begin{theorem}\label{thm:ms/conj}
$A:[a,b]^2\to [a,b]$ is continuous, associative, and is such that $a$ acts as a zero and $b$ as an identity if and only if
\begin{itemize}
\item either
$$
A(x,y) = \min(x,y) \qquad (x,y \in [a,b]),
$$
\item or there exists a continuous strictly decreasing function $f : [a,b] \to [0,+\infty]$, with $f(b) = 0$, such that
$$
A(x,y) = f^{-1}[ \, \min(f(x) + f(y),f(a)) \, ] \qquad (x,y \in [a,b]).
$$
(conjunctive Archimedean semigroup)%
\item or there exist a countable index set $K \subseteq \N$, a family of disjoint open subintervals $\{]a_k,b_k[ \, | \, k \in K\}$ of $[a,b]$
and a family $\{f_k \, | \, k \in K\}$ of continuous strictly decreasing function $f_k : [a_k,b_k] \to [0,+\infty]$, with $f_k(b_k) = 0$, such
that, for all $x,y \in [a,b]$,
$$
A(x,y) = \cases{f_k^{-1}[ \, \min(f_k(x) + f_k(y),f_k(a_k)) \, ], & if $\exists \, k \in K$ such that $x,y \in [a_k,b_k]$ \cr
        \min(x,y),                  & otherwise.\cr}
$$
(ordinal sum of conjunctive Archimedean semigroups and one-point semigroups).
\end{itemize}
\end{theorem}

One can show that ${\cal A}_{a,b,a}$ is the family of continuous, nondecreasing, and associative functions $A:[a,b]^2\to [a,b]$ such that $a$
acts as a zero and $b$ as an identity. Consequently, the description of the family ${\cal A}_{a,b,a}$ is also given by
Theorem~\ref{thm:ms/conj}. Moreover, it turns out that all functions fulfilling the assumptions of this result are symmetric, nondecreasing, and
conjunctive.

Theorem~\ref{thm:ms/conj} can also be written in a dual form as follows:

\begin{theorem}\label{thm:ms/disj}
$A:[a,b]^2\to [a,b]$ is continuous, associative, and is such that $a$ acts as an identity and $b$ as a zero if and only if
\begin{itemize}
\item either
$$
A(x,y) = \max(x,y) \qquad (x,y \in [a,b]),
$$
\item or there exists a continuous strictly increasing function $f : [a,b] \to [0,+\infty]$, with $f(a) = 0$, such that
$$
A(x,y) = f^{-1}[ \, \min(f(x) + f(y),f(b)) \, ] \qquad (x,y \in [a,b]).
$$
(disjunctive Archimedean semigroup)%
\item or there exist a countable index set $K \subseteq \N$, a family of disjoint open subintervals $\{]a_k,b_k[ \, | \, k \in K\}$ of $[a,b]$
and a family $\{f_k \, | \, k \in K\}$ of continuous strictly increasing function $f_k : [a_k,b_k] \to [0,+\infty]$, with $f_k(a_k) = 0$, such
that, for all $x,y \in [a,b]$,
$$
A(x,y) = \cases{f_k^{-1}[ \, \min(f_k(x) + f_k(y),f_k(b_k)) \, ], & if $\exists \, k \in K$ such that $x,y \in [a_k,b_k]$ \cr
        \max(x,y),                  & otherwise.\cr}
$$
(ordinal sum of disjunctive Archimedean semigroups and one-point semigroups).
\end{itemize}
\end{theorem}

As above, we can see that ${\cal A}_{a,b,b}$ is the family of continuous, nondecreasing, and associative functions $A:[a,b]^2\to [a,b]$ such
that $a$ acts as an identity and $b$ as a zero. The description of the family ${\cal A}_{a,b,b}$ is thus given by Theorem~\ref{thm:ms/disj}.
Moreover, all functions fulfilling the assumptions of this result are symmetric, nondecreasing, and disjunctive.

Theorems~\ref{thm:co/in/wid/a}, \ref{thm:ms/conj}, and \ref{thm:ms/disj}, taken together, give a complete description of the family of
continuous, nondecreasing, weakly idempotent, and associative functions $A:[a,b]^2\to [a,b]$. Imposing some additional conditions leads to the
following immediate corollaries:

\begin{corollary}
$A:[a,b]^2\to [a,b]$ is continuous, strictly increasing, weakly idempotent, and associative if and only if there exists a continuous strictly
increasing function $g : [a,b] \to [0,1]$, with $g(a) = 0$ and $g(b) = 1$, such that
\begin{itemize}
\item either
$$
A(x,y) = g^{-1}[g(x) \, g(y)] \qquad (x,y \in [a,b]),
$$
\item or
$$
A(x,y) = g^{-1}[g(x)+g(y)-g(x) \, g(y)] \qquad (x,y \in [a,b]).
$$
\end{itemize}
\end{corollary}

\begin{corollary}
$A:[a,b]^2\to [a,b]$ is symmetric, continuous, nondecreasing, weakly idempotent, and associative if and only if there exist $\alpha \in [a,b]$
and two functions $A_{a,\alpha,\alpha} \in {\cal A}_{a,\alpha,\alpha}$ and $A_{\alpha,b,\alpha} \in {\cal A}_{\alpha,b,\alpha}$ such that, for
all $x,y \in [a,b]$,
$$
A(x,y) = \cases{A_{a,\alpha,\alpha}(x,y), & if $x,y \in [a,\alpha]$ \cr
        A_{\alpha, b,\alpha}(x,y), & if $x,y \in [\alpha,b]$ \cr
        \alpha, & otherwise.\cr}
$$
\end{corollary}

\begin{corollary}\label{cor:co/in/wid/a/ident}
$A:[a,b]^2\to [a,b]$ is continuous, nondecreasing, weakly idempotent, associative, and has exactly one identity element in $[a,b]$ if and only
if $A \in {\cal A}_{a,b,a} \cup {\cal A}_{a,b,b}$.
\end{corollary}

\subsection{Internal associative functions}

We now investigate the case of internal associative functions, that is, associative means. As these functions are idempotent, we actually
investigate idempotent and associative functions. Although we have already observed that there are no continuous, strictly increasing,
idempotent, and associative functions, the class of continuous, nondecreasing, idempotent, and associative functions is nonempty and its
description can be deduced from Theorem~\ref{thm:co/in/wid/a}. However, Fodor~\cite{Fod96} had already obtained this description in a more
general framework, as follows.

\begin{theorem}\label{thm:co/in/id/a2}
Let $E$ be a real interval, finite or infinite. $A:E^2\to E$ is continuous, nondecreasing, idempotent, and associative if and only if there
exist $\alpha, \beta \in E$ such that
$$
A(x,y) = (\alpha \wedge x) \vee (\beta \wedge y) \vee (x \wedge y) \qquad ((x,y) \in E^2).
$$
\end{theorem}

Notice that, by distributivity of $\wedge$ and $\vee$, $A$ can be written also in the equivalent form:
$$
A(x,y) = (\beta \vee x) \wedge (\alpha \vee y) \wedge (x \vee y) \qquad ((x,y) \in E^2).
$$

For sequences of associative functions, the statement can be formulated as follows.

\begin{theorem}\label{thm:co/in/id/a}
Let $E$ be a real interval, finite or infinite. $(A^{(n)}:E^n\to E)_{n\ge 1}$ is an associative sequence of continuous, nondecreasing, and
idempotent functions if and only if there exist $\alpha, \beta \in E$ such that
$$
A^{(n)}(x) = (\alpha \wedge x_1) \vee \Bigl(\, \bigvee_{i=2}^{n-1}(\alpha \wedge \beta \wedge x_i)\Bigr) \vee (\beta \wedge x_n) \vee \Bigl(\,
\bigwedge_{i=1}^n x_i \Bigr) \qquad (x \in E^n, ~ n \in \N\setminus\{0\}).
$$
\end{theorem}

Before Fodor~\cite{Fod96}, the description of symmetric functions was obtained by Fung and Fu~\cite{FunFu75} and in a revisited way by Dubois
and Prade~\cite{DubPra84}. Now, the result can be formulated as follows.

\begin{theorem}\label{thm:fung/fu}
Let $E$ be a real interval, finite or infinite.\par\noindent $i)$ $A:E^2\to E$ is symmetric, continuous, nondecreasing, idempotent, and
associative if and only if there exists $\alpha \in E$ such that
$$
A(x,y) = \mathrm{median}(x,y,\alpha) \qquad (x,y \in E).
$$
$ii)$ $(A^{(n)}:E^n\to E)_{n\ge 1}$ is an associative sequence of symmetric, continuous, nondecreasing, and associative functions if and only if
there exists $\alpha \in E$ such that
\begin{equation}\label{eq:sy/co/in/id/a}
A^{(n)}(x) = \mathrm{median}\Bigl(\, \bigwedge_{i=1}^n x_i,\bigvee_{i=1}^n x_i,\alpha \Bigr) \qquad (x \in E^n, ~ n \in \N\setminus\{0\}).
\end{equation}
\end{theorem}

The previous three theorems show that the idempotency property is seldom consistent with associativity. For instance, the associative mean
(\ref{eq:sy/co/in/id/a}) is not very decisive since it leads to the predefined value $\alpha$ as soon as there exist $x_i \le \alpha$ and $x_j
\ge \alpha$.

Czoga\l a and Drewniak~\cite{CzoDre84} have examined the case when $A$ has an identity element $e \in E$, as follows.

\begin{theorem}
Let $E$ be a real interval, finite or infinite.\par\noindent $i)$ If $A:E^2\to E$ is nondecreasing, idempotent, associative, and has an identity
element $e \in E$, then there exists a decreasing function $g : E \to E$, with $g(e) = e$, such that, for all $x,y \in E$,
$$
A(x,y) = \cases{x \wedge y, & if $y < g(x)$,\cr
        x \vee y,   & if $y > g(x)$,\cr
        x \wedge y \mbox{ or } x \vee y,    & if $y = g(x)$.\cr}
$$
$ii)$ If $A:E^2\to E$ is continuous, nondecreasing, idempotent, associative, and has an identity element $e \in E$, then $A = \min$ or $\max$.
\end{theorem}

\subsection{t-norms, t-conorms, and uninorms}
\label{sec:ttu}

In fuzzy set theory, one of the main topics consists in defining fuzzy logical connectives which are appropriate extensions of logical
connectives {\bf AND}, {\bf OR}, and {\bf NOT} in the case when the valuation set is the unit interval $[0,1]$ rather than $\{0,1\}$.

Fuzzy connectives modelling {\bf AND} and {\bf OR} are called {\sl triangular norms} (t-norms for short) and {\sl triangular conorms}
(t-conorms) respectively; see Alsina et al.~\cite{AlsTriVal83} and Schweizer and Sklar~\cite{SchSkl83}.

\begin{definition}
 $i)$ A \textsl{t-norm} is a symmetric, nondecreasing, and associative function $T : [0,1]^2 \to [0,1]$ having $1$ as identity.

$ii)$ A \textsl{t-conorm} is a symmetric, nondecreasing, and associative function $S : [0,1]^2 \to [0,1]$ having $0$ as identity.
\end{definition}

The investigation of these functions has been made by Schweizer and Sklar~\cite{SchSkl61,SchSkl63} and Ling~\cite{Lin65}. There is now an
abundant literature on this topic; see the book by Klement et al.~\cite{KleMesPap00}.

Of course, the family of continuous t-norms is nothing else than the class ${\cal A}_{0,1,0}$, and the family of continuous t-conorms is the
class ${\cal A}_{0,1,1}$. Both families have been fully described in this section. Moreover, Corollary~\ref{cor:co/in/wid/a/ident} gives a
characterization of their union.

\begin{corollary}
$A:[0,1]^2\to [0,1]$ is continuous, nondecreasing, weakly idempotent, associative, and has exactly one identity in $[0,1]$ if and only if $A$ is
a continuous t-norm or a continuous t-conorm.
\end{corollary}

It is well known that t-norms and t-conorms are extensively used in fuzzy set theory, especially in modeling fuzzy connectives and implications
(see Weber~\cite{Web83}). Applications to practical problems require the use of, in a sense, the most appropriate t-norms or t-conorms. On this
issue, Fodor~\cite{Fod91} presented a method to construct new t-norms from t-norms.

It is worth noting that some properties of t-norms, such as associativity, do not play any essential role in preference modeling and choice
theory. Recently, some authors \cite{AlsMayTomTor93,DycPed84,ZimZys80} have investigated non-associative binary operation on $[0,1]$ in
different contexts. These operators can be viewed as a generalization of t-norms and t-conorms in the sense that both are contained in this kind
of operations. Moreover, Fodor~\cite{Fod91b} defined and investigated the concept of weak t-norms. His results were usefully applied to the
framework of fuzzy strict preference relations.

Further associative functions were recently introduced, namely {\sl t-operators} \cite{MasMayTor99} and {\sl uninorms} \cite{YagRyb96} (see also
\cite{MasMayTor02,MasMonTor01}), which proved to be useful in expert systems, neural networks, and fuzzy quantifiers theory.

\begin{definition}
 $i)$ A {\sl t-operator} is a symmetric, nondecreasing, associative function $F : [0,1]^2 \to [0,1]$, with $0$ and $1$ as idempotent
 elements, and such that the sections $x\mapsto F(x,0)$ and $x\mapsto F(x,1)$ are continuous on $[0,1]$.

$ii)$ A {\sl uninorm} is a symmetric, nondecreasing, and associative function $U : [0,1]^2 \to [0,1]$ having an identity.
\end{definition}

It is clear that a uninorm becomes a t-norm (resp.\ a t-conorm) when the identity is 1 (resp.\ 0).

We will not linger on this topic of t-norms, t-conorms, and uninorms. The interested reader can consult the book by Klement et
al.~\cite{KleMesPap00}. For more recent results, we also recommend an article on associative functions by Sander~\cite{San02}.

\section{Nonadditive integrals}

Many aggregation functions can be seen as nonadditive discrete integrals with respect to nonadditive measures. In this section we introduce
Choquet and Sugeno integrals. The reader can find more details on this topic in Chapter~18 of this volume.

\subsection{Motivations}

A significant aspect of aggregation in multicriteria decision making is the difference in the importance of criteria or attributes, usually modeled by using different weights. Since these weights must be taken into account during the aggregation phase, it is necessary to use
weighted functions, therefore giving up the symmetry property. Until recently, the most often used weighted aggregation functions were averaging
functions, such as the quasi-linear means (\ref{eq:co/sin/id/b}).

However, the weighted arithmetic means and, more generally, the quasi-linear means present some drawbacks. None of these functions are able to
model in an understandable way an interaction among attributes. Indeed, it is well known in multiattribute utility theory (MAUT) that these
functions lead to {\sl mutual preferential independence} among the attributes (see for instance Fishburn and Wakker~\cite{FisWak95}), which
expresses in some sense the independence of the attributes. Since these functions are not appropriate when interactive attributes are
considered, people usually tend to construct independent attributes, or attributes that are supposed to be so, causing some bias effect in
evaluation.

In order to have a flexible representation of complex interaction phenomena among attributes or criteria (e.g.\ positive or negative synergy
among some criteria), it is useful to substitute the weight vector for a nonadditive set function allowing the definition of a weight not only on each
criterion, but also on each subset of criteria.

For this purpose, the use of fuzzy measures have been proposed by Sugeno in 1974 \cite{Sug74} to generalize additive measures. It seems widely
accepted that additivity is not suitable as a required property of set functions in many real situations, due to the lack of additivity in many
facets of human reasoning. To be able to express human subjectivity, Sugeno proposed replacing the additivity property with a weaker one:
monotonicity. These non-additive monotonic measures are referred to as fuzzy measures.

We consider a discrete set of $n$ elements $N=\{1,\ldots,n\}$. Depending on the application, these elements could be players of a cooperative
game, criteria in a multicriteria decision problem, attributes, experts or voters in an opinion pooling problem, etc. To emphasize that $N$ has
$n$ elements, we will often write $N_n$.

\begin{definition}
A (discrete) {\sl fuzzy measure} on $N$ is a set function $\mu : 2^N \to [0,1]$ that is monotonic, that is $ \mu(S) \le \mu(T)$ whenever $S
\subseteq T$, and fulfills the boundary conditions $\mu(\varnothing)=0$ and $\mu(N)=1$.
\end{definition}

For any $S \subseteq N$, the coefficient $\mu(S)$ can be viewed as the weight, or importance, or strength of the combination $S$ for the
particular decision problem under consideration. Thus, in addition to the usual weights on criteria taken separately, weights on any combination
of criteria are also defined. Monotonicity then means that adding a new element to a combination cannot decrease its importance. We denote the set of fuzzy measures on $N$ as $\fn$.

When a fuzzy measure is available on $N$, it is interesting to have tools capable of summarizing all the values of a function to a single point,
in terms of the underlying fuzzy measure. These tools are the \textsl{fuzzy integrals}, a concept proposed by Sugeno~\cite{Sug74,Sug77}.

Fuzzy integrals are integrals of a real function with respect to a fuzzy measure, by analogy with Lebesgue integral which is defined with
respect to an ordinary (i.e., additive) measure. As the integral of a function in a sense represents its average value, a fuzzy integral can be
viewed as a particular case of averaging aggregation function.

Contrary to the weighted arithmetic means, fuzzy integrals are able to represent a certain kind of interaction among criteria, ranging from
redundancy (negative interaction) to synergy (positive interaction). For this reason they have been thoroughly studied in the context of
multicriteria decision problems \cite{Gra95,Gra96,Gra98,GraMurSug00}.

There are several classes of fuzzy integrals, among which the most representative are the Choquet and Sugeno integrals. In this section we
discuss these two integrals as aggregation functions. In particular, we present axiomatic characterizations for these integrals. The main
difference between them is that the former is suitable for the aggregation on interval scales, while the latter is designed for
aggregating values on ordinal scales.

\subsection{The Choquet integral}
\label{sec:cho}

The concept of Choquet integral was proposed in capacity theory \cite{Cho53}. Since then, it was used in various contexts such as nonadditive
utility theory \cite{Gil87,SarWak92,Sch86,Wak89}, theory of fuzzy measures and integrals \cite{deCLamMor91,Hoh82,MurSug89,MurSug91} (see also
the excellent edited book \cite{GraMurSug00}), but also finance \cite{DowWer92} and game theory \cite{DowWer94}.

Since this integral is viewed here as an $n$-variable aggregation function, we will adopt a function-like notation instead of the usual integral
form, and the integrand will be a set of $n$ real values, denoted $x=(x_1,\ldots,x_n)\in\R^n$.

\begin{definition}\label{de:ci}
Let $\mu \in \fn$. The (discrete) {\sl Choquet integral} of $x\in\R^n$ with respect to $\mu$ is defined by
$$
{\cal C}_{\mu}(x) := \sum_{i=1}^n x_{(i)} \, [\mu(A_{(i)})-\mu(A_{(i+1)})],
$$
where $(\cdot)$ is a permutation on $N$ such that $x_{(1)} \le \ldots \le x_{(n)}$. Also, $A_{(i)} = \{(i),\ldots,(n)\}$, and $A_{(n+1)} =
\varnothing$.
\end{definition}

For instance, if $x_3\leq x_1\leq x_2$, we have
$$
{\cal C}_{\mu}(x_1,x_2,x_3) = x_3 \, [\mu(\{3,1,2\})-\mu(\{1,2\})]+ x_1 \, [\mu(\{1,2\})-\mu(\{2\})]+ x_2 \, \mu(\{2\}).
$$

Thus, the Choquet integral is a linear expression, up to a reordering of the arguments. It is closely related to the Lebesgue integral, since
both coincide when the measure is additive:
$$
{\cal C}_{\mu}(x) = \sum_{i=1}^n \mu(i) \, x_i \qquad (x \in \R^n).
$$
In this sense, the Choquet integral is a generalization of the Lebesgue integral.

Let us now turn to axiomatizations of the Choquet integral. First of all, as we can see, this aggregation function fulfills a number of natural
properties. It is continuous, nondecreasing, unanimously increasing, idempotent, internal, and meaningful for the same input-output interval
scales; see for instance Grabisch~\cite{Gra96}. It also fulfills the \textsl{comonotonic additivity} property \cite{Del71,Sch86}, that is,
$$
f(x_1+x'_1,\ldots,x_n+x'_n)=f(x_1,\ldots,x_n)+f(x'_1,\ldots,x'_n)
$$
for all \textsl{comonotonic} vectors $x,x'\in\R^n$. Two vectors $x,x'\in\R^n$ are comonotonic if there exists a permutation $\sigma$ on
$N$ so that
$$
x_{\sigma(1)}\le\cdots\le x_{\sigma(n)} \quad\mbox{and}\quad x'_{\sigma(1)}\le\cdots\le x'_{\sigma(n)}.
$$
An interpretation of this property in multicriteria decision making can be found in Modave et al.~\cite{ModDubGraPra97,ModGra98}.

The following result \cite[Proposition 4.1]{MarMatTou99} gives a characterization of the 2-variable Choquet integral in a very natural way.

\begin{proposition}
$f:\R^2\to\R$ is nondecreasing and meaningful for the same input-output interval scales if and only if there exists $\mu\in {\cal F}_2$ such
that $f={\cal C}_{\mu}$.
\end{proposition}

The class of $n$-variable Choquet integrals has been first characterized by Schmeidler~\cite{Sch86} by using monotonic additivity; see also
\cite{deCJor92}, \cite{deCLamMor91}, \cite{Gra93}, and \cite[Theorem 8.6]{GraNguWal95}. Note that this result was stated and proved in the
continuous case (infinite) instead of the discrete case.

\begin{theorem}
$f:\R^n\to\R$ is nondecreasing, comonotonic additive, and fulfills $f(\mathbf{1}_N)=1$ if and only if there exists $\mu\in \fn$ such that $f={\cal
C}_\mu$.
\end{theorem}

Since the Choquet integral is defined from a fuzzy measure, it is sometimes useful to consider, for a given set $N$, the family of Choquet
integrals on $N$ as a set of functions
$$
\{f_{\mu}:\R^n\to\R\mid \mu\in\fn\}
$$
or, equivalently, as a function $f:\R^n\times\fn\to\R$.

Let us mention a first characterization of the family of Choquet integrals on $N$; see Groes et al.~\cite{GroJacSloTra98}. For any $S\subseteq
N$, $S\neq\varnothing$, denote by $\mu_S$ the fuzzy measure on $N$ defined by $\mu_S(T)=1$ if $T\supseteq S$ and 0 otherwise.

\begin{theorem}
The class of functions $\{f_{\mu}:\R^n \to \R \mid \mu\in \fn\}$ fulfills the following properties:
\begin{itemize}
\item for any $\mu,\nu\in\fn$ and any $\lambda\in\R$ such that $\lambda\mu+(1-\lambda)\nu\in\fn$, we have
$$f_{\lambda \mu+(1-\lambda)\nu}=\lambda f_{\mu}+(1-\lambda)f_{\nu}\, ,$$
\item for any $S\subseteq N$, we have $f_{\mu_S}=\min_S$,
\end{itemize}
if and only if $f_{\mu}={\cal C}_\mu$ for all $\mu\in\fn$.
\end{theorem}

A second characterization obtained by the author \cite{Mar98,Mar00g}, can be stated as follows.

\begin{theorem}
The class of functions $\{f_{\mu}:\R^n \to \R \mid \mu\in \fn\}$ fulfills the following properties:
\begin{itemize}
\item any function $f_{\mu}$ is a linear expression of $\mu$, that is, there exist $2^n$ functions $g_T:\R^n\to \R$ $(T\subseteq N)$ such that
$f_{\mu}=\sum_{T\subseteq N} g_T \, \mu(T)$ for all $\mu\in\fn$,

\item for any $\mu\in\fn$ and any $S\subseteq N$, we have $f_{\mu}(\mathbf{1}_S)=\mu(S),$

\item for any $\mu\in\fn$, the function $f_{\mu}$ is nondecreasing and meaningful for the same input-output interval scales,
\end{itemize}
if and only if $f_{\mu}={\cal C}_{\mu}$ for all $\mu\in\fn$.
\end{theorem}

These two characterizations are natural and similar to each other. The linearity condition proposed in the second
characterization is useful if we want to keep the aggregation model as simple as possible. Technically, this condition is equivalent to the
superposition condition, that is,
$$
f_{\lambda_1\mu+\lambda_2\nu} = \lambda_1 f_{\mu}+\lambda_2 f_{\nu}
$$
for all $\mu,\nu\in\fn$ and all $\lambda_1,\lambda_2\in\R$ such that $\lambda_1\mu+\lambda_2\nu\in\fn$. Of course, linearity implies the first
condition of the first characterization. Moreover, under this linearity condition, the other conditions are equivalent. In fact, in the proof of
the second characterization \cite{Mar98,Mar00g}, the author replaced the condition $f_{\mu_S}=\min_S$ with the three conditions:
$f_{\mu}(\mathbf{1}_S)=\mu(S)$, nondecreasing monotonicity, and meaningfulness for the same input-output interval scales of $f_\mu$.

We also have the following three results \cite[\S4.2.3]{Mar98}:

\begin{proposition}
The Choquet integral ${\cal C}_{\mu}:\R^n\to\R$ is bisymmetric if and only if
$$
{\cal C}_{\mu}\in \{{\rm min}_S,{\rm max}_S\mid S\subseteq N\}\cup \{{\rm WAM}_{\omega}\mid \omega\in [0,1]^n\}.
$$
\end{proposition}

\begin{proposition}
The sequence of Choquet integrals ${\cal C}:=({\cal C}^{(n)}_{\mu^{(n)}}:\R^n\to\R)_{n\ge 1}$ is bisymmetric if and only if
\begin{itemize}
\item either, for any $n\in\N\setminus\{0\}$, there exists $S\subseteq N_n$ such that ${\cal C}^{(n)}_{\mu^{(n)}}={\rm min}_S$,

\item or, for any $n\in\N\setminus\{0\}$, there exists $S\subseteq N_n$ such that ${\cal C}^{(n)}_{\mu^{(n)}}={\rm max}_S$,

\item or, for any $n\in\N\setminus\{0\}$, there exists $\omega\in [0,1]^n$ such that ${\cal C}^{(n)}_{\mu^{(n)}}={\rm WAM}_{\omega}$.
\end{itemize}
\end{proposition}

\begin{proposition}
The sequence of Choquet integrals ${\cal C}:=({\cal C}^{(n)}_{\mu^{(n)}}:\R^n\to\R)_{n\ge 1}$ is decomposable if and only if
\begin{itemize}
\item either ${\cal C}=({\rm min}^{(n)})_{n\ge 1}$,

\item or ${\cal C}=({\rm max}^{(n)})_{n\ge 1}$,

\item or there exists $\theta\in [0,1]$ such that, for any $n\in\N\setminus\{0\}$, we have ${\cal C}^{(n)}_{\mu^{(n)}}={\rm WAM}_{\omega}$, with
$$
\omega_i = \frac{(1-\theta)^{n-i}\theta^{i-1}}{\sum_{j=1}^n (1-\theta)^{n-j}\theta^{j-1}}\qquad (i\in N_n).
$$
\end{itemize}
\end{proposition}

\begin{proposition}
The sequence of Choquet integrals ${\cal C}:=({\cal C}^{(n)}_{\mu^{(n)}}:\R^n\to\R)_{n\ge 1}$ is associative if and only if
$$
{\cal C}=({\rm min}^{(n)})_{n\ge 1} \mbox{ or }({\rm max}^{(n)})_{n\ge 1}\mbox{ or } (P_1^{(n)})_{n\ge 1}\mbox{ or } (P_n^{(n)})_{n\ge 1}.
$$
\end{proposition}

Let us now consider certain special cases of the Choquet integral, namely the weighted arithmetic means (WAM) and the ordered weighted
averaging functions (OWA).

The {\sl weighted arithmetic mean}\/ ${\rm WAM}_{\omega}$ is a Choquet integral defined from an additive measure. It fulfills the classical
\textsl{additivity} property:
$$
f(x_1+x'_1,\ldots,x_n+x'_n)=f(x_1,\ldots,x_n)+f(x'_1,\ldots,x'_n)
$$
for all vectors $x,x'\in\R^n$. More exactly, we have the following results (see Marichal~\cite[\S4.2.4]{Mar98} and Murofushi and
Sugeno~\cite{MurSug93}).

\begin{proposition}
The Choquet integral ${\cal C}_{\mu}:\R^n\to\R$ is additive if and only if there exists $\omega\in [0,1]^n$ such that ${\cal C}_{\mu}={\rm
WAM}_{\omega}$.
\end{proposition}

\begin{proposition}
$A:\R^n\to\R$ is nondecreasing, meaningful for the same input-output interval scales, and additive if and only if there exists $\omega\in
[0,1]^n$ such that $A={\rm WAM}_{\omega}$.
\end{proposition}

The {\sl ordered weighted averaging functions}\/ ${\rm OWA}_{\omega}$ were proposed in 1988 by Yager~\cite{Yag88}. Since their introduction,
these aggregation functions have been applied to many fields as neural networks, data base systems, fuzzy logic controllers, and group decision
making. An overview on these functions can be found in the book edited by Yager and Kacprzyk~\cite{YagKac97}; see also Grabisch et
al.~\cite{GraMurSug00}.

The following result, ascribed to Grabisch~\cite{Gra95b} (see \cite{Mar02b} for a concise proof), shows that the OWA function is nothing other
than a Choquet integral with respect to a {\sl cardinality-based} fuzzy measure i.e.\ a fuzzy measure depending only on the cardinalities of
the subsets.

\begin{proposition}\label{prop:owa}
Let $\mu\in\fn$. The following assertions are equivalent:\par
\begin{tabular}{rl}
  i) & For any $S, S' \subseteq N$ such that $|S| = |S'|$, we have $\mu(S) = \mu(S')$.\\
 ii) & There exists a weight vector $\omega$ such that ${\cal C}_{\mu} = {\rm OWA}_{\omega}$. \\
iii) & ${\cal C}_{\mu}$ is a symmetric function.
\end{tabular}
\end{proposition}

The fuzzy measure $\mu$ associated to an ${\rm OWA}_\omega$ is given by
$$
\mu(S) = \sum_{i=n-s+1}^n \omega_i\qquad (S \subseteq N, \enskip S \neq \varnothing).
$$
Conversely, the weights associated to ${\rm OWA}_\omega$ are given by
$$
\omega_{n-s} = \mu(S\cup i)-\mu(S) \qquad (i \in N, \enskip S \subseteq N\setminus i).
$$

The class of OWA functions includes an important subfamily, namely the order statistics
$$
{\rm OS}_k(x) = x_{(k)},
$$
when $\omega_k = 1$ for some $k \in N$. In this case, we have, for any $S\subseteq N$,
$$
\mu(S) = \cases{ 1, & if $s\ge n-k+1$, \cr 0, & otherwise.\cr }
$$
This subfamily itself contains the minimum, the maximum, and the median.

Axiomatizations of the class of OWA functions can be immediately derived from those of the Choquet integral and from Proposition~\ref{prop:owa}.

\subsection{The Sugeno integral}

The Sugeno integral \cite{Sug74,Sug77} was introduced as a fuzzy integral, that is, an integral defined from a fuzzy measure. This integral has
then been thoroughly investigated and used in many domains (an overview can be found in Dubois et al.~\cite{DubMarPraRouSab01} and the volume edited
by Grabisch et al.~\cite{GraMurSug00}).

As for the Choquet integral, we give here the definition of the Sugeno integral in its discrete (finite) version, which is nothing other
than an aggregation function from  $[0,1]^n$ into $[0,1]$.

\begin{definition}\label{de:si}
Let $\mu \in \fn$. The (discrete) {\sl Sugeno integral} of $x\in [0,1]^n$ with respect to $\mu$ is defined by
$$
{\cal S}_{\mu}(x) := \bigvee_{i=1}^n [x_{(i)} \wedge \mu(A_{(i)})],
$$
where $(\cdot)$ is a permutation on $N$ such that $x_{(1)} \le \ldots \le x_{(n)}$. Also, $A_{(i)} = \{(i),\ldots,(n)\}$, and $A_{(n+1)} =
\varnothing$.
\end{definition}

Exactly as in the definition of the Choquet integral, the ``coefficient'' associated with each variable $x_i$ is fixed uniquely by the
permutation $(\cdot)$. For instance, if $x_3\leq x_1\leq x_2$, then we have
$$
{\cal S}_{\mu}(x_1,x_2,x_3) = [x_3 \wedge \mu(\{3,1,2\})] \vee [x_1 \wedge \mu(\{1,2\})] \vee [x_2 \wedge \mu(\{2\})].
$$

From the definition, we can immediately deduce that
$$
{\cal S}_{\mu}(x) \in \{x_1,\ldots,x_n\} \cup \{\mu(S) \, | \, S \subseteq N\} \qquad (x \in [0,1]^n).
$$
Moreover, similarly to the Choquet integral, we have
$$
{\cal S}_{\mu}(\mathbf{1}_S) = \mu(S) \qquad (S \subseteq N),
$$
which shows the athe Sugeno integral is completely determined par its values at the vertices of the hypercube $[0,1]^n$.

It was proved \cite{Gre87,Mar00c,Sug74} that the Sugeno integral can also be set in the following form, which does not require the reordering of
the variables:
$$
{\cal S}_{\mu}(x) = \bigvee_{T \subseteq N} \Big[\mu(T) \wedge (\bigwedge_{i \in T} x_i)\Big] \qquad (x \in [0,1]^n).
$$

It was also proved \cite{KanBya78} that the Sugeno integral is a kind of weighted median:
$$
{\cal S}_{\mu}(x) = \mathrm{median}[x_1,\ldots,x_n,\mu(A_{(2)}),\mu(A_{(3)}),\ldots,\mu(A_{(n)})] \qquad (x \in [0,1]^n).
$$
For instance, if $x_3\le x_1\le x_2$, then
$$
{\cal S}_{\mu}(x_1,x_2,x_3)=\mathrm{median}[x_1,x_2,x_3,\mu(1,2),\mu(2)].
$$

The following result \cite{Mar01} shows that the Sugeno integral is a rather natural concept and, contrary to the Choquet integral, it is
suitable for an aggregation in an ordinal context.

\begin{proposition}
Any weakly idempotent function $A:[0,1]^n\to [0,1]$, whose (well-formed) expression is made up of variables $x_1,\ldots,x_n$, constants
$r_1,\ldots,r_m\in [0,1]$, lattice operations $\wedge=\min$ et $\vee=\max$, and parentheses is a Sugeno integral (and conversely).
\end{proposition}

Let us now turn to axiomatizations of the Sugeno integral. We can easily see that the Sugeno integral is a continuous, nondecreasing,
unanimously increasing, idempotent, and internal function. It also fulfills the \textsl{comonotonic minitivity} and \textsl{comonotonic
maxitivity} properties \cite{deCLamMor91}, that is
\begin{eqnarray*}
f(x_1\wedge x'_1,\ldots,x_n \wedge x'_n) &=& f(x_1,\ldots,x_n)\wedge f(x'_1,\ldots,x'_n)\\
f(x_1\vee x'_1,\ldots,x_n \vee x'_n) &=& f(x_1,\ldots,x_n)\vee f(x'_1,\ldots,x'_n)
\end{eqnarray*}
for all comonotonic vectors $x,x'\in [0,1]^n$. More specifically, it is \textsl{weakly minitive} and \textsl{weakly maxitive}, that is, it
fulfills
\begin{eqnarray*}
f(x_1\wedge r,\ldots,x_n \wedge r) &=& f(x_1,\ldots,x_n)\wedge r\\
f(x_1\vee r,\ldots,x_n \vee r) &=& f(x_1,\ldots,x_n)\vee r
\end{eqnarray*}
for all vectors $x\in [0,1]^n$ and all $r\in [0,1]$. Even more specifically, by replacing $x$ with the Boolean vector $\mathbf{1}_S$ in these
two equations above, one see that it is also {\sl non-compensative}, that is, it fulfills
$$
f(r\mathbf{1}_S)\in \{f(\mathbf{1}_S),r\} ~~\mbox{and}~~ f(\mathbf{1}_S+r\mathbf{1}_{N\setminus S})\in \{f(\mathbf{1}_S),r\}
$$
for all $S\subseteq N$ and all $r\in [0,1]$.

Comonotonic minitivity and maxitivity have been interpreted in the context of aggregation of fuzzy subsets by Ralescu and Ralescu
\cite{RalRal97}. Non-compensation has been interpreted in decision making under uncertainty in Dubois et al.~\cite{DubMarPraRouSab01}.

The main axiomatizations of the Sugeno integral as an aggregation function are summarized in the following result; see
Marichal~\cite{Mar98,Mar00c}:

\begin{theorem}
Let $A:[0,1]^n\to [0,1]$. The following assertions are equivalent:
\begin{itemize}
\item $A$ is nondecreasing, idempotent, and non-compensative,

\item $A$ is nondecreasing, weakly minitive and weakly maxitive,

\item $A$ is nondecreasing, idempotent, comonotonic minitive and maxitive,

\item There exists $\mu\in\fn$ such that $A={\cal S}_{\mu}$.
\end{itemize}
\end{theorem}

The 2-variable Sugeno integral can be characterized in a very natural way by means of the associativity property. Indeed,
Theorem~\ref{thm:co/in/id/a2} can be rewritten as follows.

\begin{proposition}
$A:[0,1]^2\to [0,1]$ is continuous, nondecreasing, idempotent, and associative if and only if there exists $\mu\in {\cal F}_2$ such that
$A={\cal S}_{\mu}$.
\end{proposition}

Considering associative or decomposable sequences, we have the following result; see Marichal~\cite[p.~113]{Mar98}.

\begin{proposition}
Let $A:=(A^{(n)}:[0,1]^n\to [0.1])_{n\ge 1}$ be a sequence of functions. Then the following assertions are equivalent:
\begin{itemize}
\item $A$ is an associative sequence of Sugeno integrals.

\item $A$ is a decomposable sequence of Sugeno integrals.

\item $A$ is an associative sequence of continuous, nondecreasing, and idempotent functions.

\item There exist $\alpha, \beta \in [0,1]$ such that
$$
A^{(n)}(x) = (\alpha \wedge x_1) \vee \Bigl(\, \bigvee_{i=2}^{n-1}(\alpha \wedge \beta \wedge x_i)\Bigr) \vee (\beta \wedge x_n) \vee \Bigl(\,
\bigwedge_{i=1}^n x_i \Bigr) \qquad (x \in [0,1]^n, ~ n \in \N\setminus\{0\}).
$$
\end{itemize}
\end{proposition}

Just as the Choquet integral includes two main subclasses, namely the weighted arithmetic means and the ordered weighted averaging functions,
the Sugeno integral includes the weighted minimum and maximum and the ordered weighted minimum and maximum. These functions have been introduced
and investigated in Dubois and Prade~\cite{DubPra86} and Dubois et al.~\cite{DubPraTes88}, respectively.

For any vector $\omega =(\omega_1,\ldots,\omega_n)\in [0,1]^n$ such that $\bigvee_{i=1}^n \omega_i=1$, the {\sl weighted maximum} associated
with $\omega$ is defined by
$$
{\rm pmax}_{\omega}(x)=\bigvee_{i=1}^n(\omega_i\wedge x_i) \qquad (x\in [0,1]^n).
$$
For any vector $\omega=(\omega_1,\ldots,\omega_n)\in [0,1]^n$ such that $\bigwedge_{i=1}^n \omega_i=0$, the {\sl weighted minimum} associated
with $\omega$ is defined by
$$
{\rm pmin}_{\omega}(x)=\bigwedge_{i=1}^n(\omega_i\vee x_i) \qquad (x\in [0,1]^n).
$$

The functions ${\rm pmax}_{\omega}$ and ${\rm pmin}_{\omega}$ can be characterized as follows; see \cite{DubPra86,Mar98,RalSug96}:

\begin{proposition}
Let $\mu\in\fn$. The following assertions are equivalent:
\begin{itemize}
\item $\mu$ is a possibility measure, that is such that
$$
\mu(S\cup T)=\mu(S)\vee\mu(T)\qquad (S,T\subseteq N).
$$

\item There exists $\omega\in [0,1]^n$ such that ${\cal S}_{\mu}={\rm pmax}_{\omega}$.

\item ${\cal S}_{\mu}(x_1\vee x'_1,\ldots,x_n\vee x'_n)={\cal S}_{\mu}(x_1,\ldots,x_n)\vee {\cal S}_{\mu}(x'_1,\ldots,x'_n) \qquad (x,x'\in [0,1]^n)$.
\end{itemize}
\end{proposition}

\begin{proposition}
Let $\mu\in\fn$. The following assertions are equivalent:
\begin{itemize}
\item $\mu$ is a necessity measure, that is such that
$$
\mu(S\cap T)=\mu(S)\wedge\mu(T)\qquad (S,T\subseteq N).
$$
\item There exists $\omega\in [0,1]^n$ such that ${\cal S}_{\mu}={\rm pmin}_{\omega}$.

\item ${\cal S}_{\mu}(x_1\wedge x'_1,\ldots,x_n\wedge x'_n)={\cal S}_{\mu}(x_1,\ldots,x_n)\wedge {\cal S}_{\mu}(x'_1,\ldots,x'_n) \qquad (x,x'\in
[0,1]^n)$.
\end{itemize}
\end{proposition}

For any vector $\omega =(\omega_1,\ldots,\omega_n)\in [0,1]^n$ such that $\bigvee_{i=1}^n \omega_i=1$, the {\sl ordered weighted maximum}
associated with $\omega$ is defined by
$$
{\rm opmax}_{\omega}(x)=\bigvee_{i=1}^n(\omega_i\wedge x_{(i)}) \qquad (x\in [0,1]^n).
$$
For any vector $\omega=(\omega_1,\ldots,\omega_n)\in [0,1]^n$ such that $\bigwedge_{i=1}^n \omega_i=0$, the {\sl ordered weighted minimum}
associated with $\omega$ is defined by
$$
{\rm opmin}_{\omega}(x)=\bigwedge_{i=1}^n(\omega_i\vee x_{(i)}) \qquad (x\in [0,1]^n).
$$

Surprisingly enough, the class of ordered weighted minima coincides with that of ordered weighted maxima and identifies with the symmetric
Sugeno integrals. The result can be stated as follows; see \cite{Gra95b,Mar98}.

\begin{proposition}
Let $\mu\in\fn$. The following assertions are equivalent:
\begin{itemize}
\item $\mu$ depends only of the cardinalities of the subsets.

\item There exists $\omega\in [0,1]^n$ such that ${\cal S}_{\mu}={\rm opmax}_{\omega}$.

\item There exists $\omega\in [0,1]^n$ such that ${\cal S}_{\mu}={\rm opmin}_{\omega}$.

\item ${\cal S}_{\mu}$ is a symmetric function.
\end{itemize}
\end{proposition}

Using the fact that the Sugeno integral is also a weighted median, we can write
\begin{eqnarray*}
{\rm opmax}_{\omega}(x) &=& \mathrm{median}(x_1,\ldots,x_n,\omega_2,\ldots,\omega_n),\\
{\rm opmin}_{\omega}(x) &=& \mathrm{median}(x_1,\ldots,x_n,\omega_1,\ldots,\omega_{n-1}).
\end{eqnarray*}

Another interesting subclass is that of lattice polynomials, which are nothing other than Sugeno integrals defined from fuzzy measures taking
their values in $\{0,1\}$. We will characterize these functions in the final section.

\section{Aggregation on ratio and interval scales}

In this section, we present the families of aggregation functions that are meaningful for ratio and interval scales (see
Definition~\ref{de:signi}). First of all, we have the following two results concerning ratio scales; see \cite[Chapter 20]{AczDho89}, \cite[p.~439]{AczGroSch94}, and
\cite[case~\#2]{AczRobRos86}.
\begin{theorem}
$A:\,]0,\infty[^n\to\, ]0,\infty[$ is meaningful for the same input-output ratio scales if and only if
$$
A(x)=x_1\, F\Big(\frac{x_2}{x_1},\ldots,\frac{x_n}{x_1}\Big) \qquad (x\in ]0,\infty[^n),
$$
with $F:\,]0,\infty[^{n-1}\to\, ]0,\infty[$ arbitrary ($F=$ constant if $n=1$).
\end{theorem}

\begin{theorem}
$A:\,]0,\infty[^n\to\, ]0,\infty[$ is meaningful for the same input ratio scales if and only if
$$
A(x)=g(x_1)\, F\Big(\frac{x_2}{x_1},\ldots,\frac{x_n}{x_1}\Big) \qquad (x\in\, ]0,\infty[^n),
$$
with $F:\, ]0,\infty[^{n-1}\to\, ]0,\infty[$ arbitrary ($F=$ constant if $n=1$) and $g:\,]0,\infty[\to\, ]0,\infty[$ such that $g(xy)=g(x)g(y)$
for all $x,y\in\, ]0,\infty[$. $g(x)=x^c$ if $A$ is continuous ($c$ arbitrary).
\end{theorem}

We have the following results regarding interval scales; see \cite[case~\#5]{AczRobRos86} and \cite[\S3.4.1]{Mar98}.

\begin{theorem}
$A:\R^n\to \R$ is meaningful for the same input-output interval scales if and only if
$$
A(x) = \cases{ {\rm S}(x) \, F\Bigl(\frac{x_1-{\rm AM}(x)}{{\rm S}(x)},\ldots,\frac{x_n-{\rm AM}(x)}{{\rm S}(x)}\Bigr) + {\rm AM}(x), & if ${\rm
S}(x) \neq 0$,\cr x_1, & if ${\rm S}(x) = 0$,\cr}
$$
where ${\rm S}(x) = \sqrt{\sum_{i=1}^n (x_i - {\rm AM}(x))^2}$ and $F:\R^n\to \R$ arbitrary ($A(x) = x$ if $n=1$).
\end{theorem}

\begin{theorem}
$A:\R^n\to \R$ is meaningful for the same input interval scales if and only if
$$
A(x) = \cases{ {\rm S}(x) \, F\Bigl(\frac{x_1-{\rm AM}(x)}{{\rm S}(x)},\ldots,\frac{x_n-{\rm AM}(x)}{{\rm S}(x)}\Bigr) + a \, {\rm AM}(x) + b, &
if ${\rm S}(x) \neq 0$,\cr a\, x_1 + b, & if ${\rm S}(x) = 0$,\cr}
$$
or
$$
A(x) = \cases{ g({\rm S}(x)) \, F\Bigl(\frac{x_1-{\rm AM}(x)}{{\rm S}(x)},\ldots,\frac{x_n-{\rm AM}(x)}{{\rm S}(x)}\Bigr) + b, & if ${\rm S}(x)
\neq 0$,\cr b, & if ${\rm S}(x) = 0$,\cr}
$$
where $a, b \in \R$, ${\rm S}(x) = \sqrt{\sum_{i=1}^n (x_i - {\rm AM}(x))^2}$, $F:\R^n\to \R$ arbitrary ($A(x) = a x + b$ if $n=1$), and
$g:\R\to\, ]0,\infty[$ such that $g(x y)=g(x)g(y)$ for all $x, y \in \R$.
\end{theorem}

The restriction of these families to nondecreasing functions and strictly increasing functions is discussed in Acz\'el et
al.~\cite{AczGroSch94}.

In the rest of this section, we present axiomatizations of some subfamilies of functions that are meaningful for the same input-output interval
scales (these results are extracted from Marichal et al.~\cite{MarMatTou99}). For instance, we observed in Section~\ref{sec:cho} that
the discrete Choquet integral fulfills this property. More generally, it is clear that any aggregation function obtained by composition of an
arbitrary number of discrete Choquet integrals is again meaningful for the same input-output interval scales. These functions, called {\sl
composed Choquet integrals}, have been investigated for instance in Narukawa and Murofushi~\cite{NarMur02}.

If we confine ourselves to bisymmetric functions, we have the following results.

\begin{proposition}
$A:\R^n\to\R$ is nondecreasing, meaningful for the same input-output interval scales, and bisymmetric if and only if
$$
A\in\{{\rm min}_S, {\rm max}_S\mid S\subseteq N\}\cup\{{\rm WAM}_{\omega}\mid \omega\in [0,1]^n\}.
$$
\end{proposition}

\begin{corollary}
$A:\R^n\to\R$ is symmetric, nondecreasing, meaningful for the same input-output interval scales, and bisymmetric if and only if
$$A\in\{{\rm min}, {\rm max}, {\rm AM}\}.$$
\end{corollary}

\begin{proposition}
$(A^{(n)}:\R^n\to\R)_{n\ge 1}$ is a bisymmetric sequence of nondecreasing and meaningful functions for the same input-output interval scales if
and only if
\begin{itemize}
\item either, for any $n\in\N\setminus\{0\}$, there exists $S\subseteq N_n$ such that $M^{(n)}={\rm min}_S$,

\item or, for any $n\in\N\setminus\{0\}$, there exists $S\subseteq N_n$ such that $M^{(n)}={\rm max}_S$,

\item or, for any $n\in\N\setminus\{0\}$, there exists $\omega\in [0,1]^n$ such that $M^{(n)}={\rm WAM}_{\omega}$.
\end{itemize}
\end{proposition}

\begin{corollary}
$A:=(A^{(n)}:\R^n\to\R)_{n\ge 1}$ is a bisymmetric sequence of symmetric, nondecreasing, and meaningful functions for the same input-output
interval scales if and only if
$$
A = ({\rm min}^{(n)})_{n\ge 1} \mbox{ or } ({\rm max}^{(n)})_{n\ge 1} \mbox{ or } ({\rm AM}^{(n)})_{n\ge 1}.
$$
\end{corollary}

Let us now consider the decomposable and associative sequences of aggregation functions. We have the following results.

\begin{proposition}
$A:=(A^{(n)}:\R^n\to\R)_{n\ge 1}$ is a decomposable sequence of nondecreasing and meaningful functions for the same input-output interval scales
if and only if
\begin{itemize}
\item either $A=({\rm min}^{(n)})_{n\ge 1}$,

\item or $A=({\rm max}^{(n)})_{n\ge 1}$,

\item or there exists $\theta\in [0,1]$ such that, for any $n\in\N\setminus\{0\}$, we have $A^{(n)}={\rm WAM}_{\omega}$, with
$$
\omega_i =\frac{(1-\theta)^{n-i}\theta^{i-1}}{\sum_{j=1}^n (1-\theta)^{n-j}\theta^{j-1}}\qquad (i\in N_n).
$$
\end{itemize}
\end{proposition}

\begin{corollary}
$A:=(A^{(n)}:\R^n\to\R)_{n\ge 1}$ is a decomposable sequence of symmetric, nondecreasing, and meaningful functions for the same input-output
interval scales if and only if
$$
A = ({\rm min}^{(n)})_{n\ge 1} \mbox{ or } ({\rm max}^{(n)})_{n\ge 1} \mbox{ or } ({\rm AM}^{(n)})_{n\ge 1}.
$$
\end{corollary}

\begin{proposition}
$A:=(A^{(n)}:\R^n\to\R)_{n\ge 1}$ is an associative sequence of nondecreasing and meaningful functions for the same input-output interval scales
if and only if
$$
A = ({\rm min}^{(n)})_{n\ge 1} \mbox{ or } ({\rm max}^{(n)})_{n\ge 1} \mbox{ or } ({\rm P}_1^{(n)})_{n\ge 1} \mbox{ or } ({\rm P}_n^{(n)})_{n\ge
1}.
$$
\end{proposition}

\begin{corollary}
$A:=(A^{(n)}:\R^n\to\R)_{n\ge 1}$ is an associative sequence of symmetric, nondecreasing, and meaningful functions for the same input-output
interval scales if and only if
$$
A = ({\rm min}^{(n)})_{n\ge 1} \mbox{ or } ({\rm max}^{(n)})_{n\ge 1}.
$$
\end{corollary}

\section{Aggregation on ordinal scales}
\label{sec:agrord}

In this final section, we consider aggregation functions that are meaningful for the same input-output ordinal scales. Their description is not
immediate and requires the concept of invariant sets. Denote the set of strictly increasing bijections of $\R$ by $\Phi$.

\begin{definition}
A nonempty subset $I\subseteq\R^n$ is said to be {\sl invariant} if
$$
x\in I ~~\Rightarrow ~~\phi(x)\in I \qquad (\phi\in\Phi).
$$
Such a set is said to be {\sl minimal} if it does not contain any proper invariant subset.
\end{definition}

The family ${\cal I}$ of all invariant subsets of $\R^n$ provides a partition of $\R^n$ into equivalence classes, where $x,y\in\R^n$ are
equivalent if there exists $\phi\in\Phi$ such that $y=\phi(x)$. In fact, one can show that any invariant subset is of the form
$$
I=\{x\in\R^n\mid x_{\pi(1)} \vartriangleleft_1 \cdots \vartriangleleft_{n-1} x_{\pi(n)}\},
$$
where $\pi\in\Pi_N$ and $\vartriangleleft_i\in\{<,\le\}$ for $i=1,\ldots,n-1$.

The meaningful functions for the same input-output ordinal scales have been investigated by many authors \cite{Mar02c,MarRou93,MesRuc,Ovc9798}.
They can be described as follows \cite{Ovc9798}.

\begin{theorem}
$A:\R^n\to\R$ is meaningful for the same input-output ordinal scales if and only if, for any $I\in {\cal I}$, there exists $i\in N$ such that
$A|_I=P_i|_I$ is the $i$th projection.
\end{theorem}

The meaningful functions for the same input ordinal scales have also been widely studied \cite{Mar02c,MarMesRuc,Orl81,Ovc96,Yan89}. They have
been described as follows \cite{MarMesRuc}.

\begin{theorem}
$A:\R^n\to\R$ is meaningful for the same input ordinal scales if and only if, for any $I\in {\cal I}$, there exists $i_I\in N$ and a constant or
strictly monotone function $g_I:P_{i_I}(I)\to\R$ such that
$$
A|_I=g_I\circ P_{i_I}|_I,
$$
where, for any $I,J\in {\cal I}$, either $g_I=g_J$, or ${\rm ran}(g_I)={\rm ran}(g_J)$ est un singleton, or ${\rm ran}(g_I)<{\rm ran}(g_J)$, or
${\rm ran}(g_I)>{\rm ran}(g_J)$.
\end{theorem}

We therefore see that the meaningful functions for the same input-output ordinal scales reduce to projections on each invariant subset. In addition, the
meaningful functions for the same input ordinal scales reduce to constants or transformed projections on these invariant subsets.

The restriction of these functions to nondecreasing and/or continuous functions has also been studied. To describe these subfamilies, we need
the concept of lattice polynomials.

\begin{definition}
A {\sl lattice polynomial} of $n$ variables is a well-formed expression involving $n$ variables $x_1,\ldots,x_n$ linked by the lattice
operations $\wedge=\min$ and $\vee=\max$ in an arbitrary combination of parentheses.
\end{definition}

For instance, $L(x)=(x_1\vee x_2)\wedge x_3$ is a 3-variable lattice polynomial.

One can show (see Birkhoff~\cite[Chapter 2, \S5]{Bir67}) that any $n$-variable lattice polynomial can be written in disjunctive form as
$$
L_{\gamma}(x)=\bigvee_{\textstyle{S\subseteq N \atop \gamma(S)=1}}\bigwedge_{i\in S}x_i \qquad (x\in\R^n),
$$
where $\gamma:2^N\to \{0,1\}$ is a binary fuzzy measure (i.e.\ with values in $\{0,1\}$). We denote the family of these fuzzy
measures on $N$ by $\Gamma_N$.

It was also proved \cite{Mar01} that the class of lattice polynomials restricted to the domain $[0,1]^n$ identifies with the intersection
between the family of Choquet integrals on $[0,1]^n$ and the family of Sugeno integrals.

Regarding nondecreasing functions, we have the following descriptions \cite{Mar02c,MarMesRuc}.

\begin{proposition}
$A:\R^n\to\R$ is nondecreasing and meaningful for the same input-output ordinal scales if and only if there exists $\gamma\in\Gamma_N$ such that
$A=L_{\gamma}$.
\end{proposition}

\begin{proposition}
$A:\R^n\to\R$ is nondecreasing and meaningful for the same input ordinal scales if and only if there exists $\gamma\in\Gamma_N$ and a constant
or strictly monotone function $g:\R\to\R$ such that $A=g\circ L_{\gamma}$.
\end{proposition}

The functions in the above two theorems are continuous, up to discontinuities of function $g$.

Regarding continuous functions, we have the following results \cite{Mar02c}.

\begin{corollary}
$A:\R^n\to\R$ is continuous and meaningful for the same input-output ordinal scales if and only if there exists $\gamma\in\Gamma_N$ such that
$A=L_{\gamma}$.
\end{corollary}

\begin{corollary}
$A:\R^n\to\R$ is continuous and meaningful for the same input ordinal scales if and only if there exists $\gamma\in\Gamma_N$ and a constant or
continuous and strictly monotone function $g:\R\to\R$ such that $A=g\circ L_{\gamma}$.
\end{corollary}

Lattice polynomials are idempotent, but not necessarily symmetric. Actually, symmetric lattice polynomials are exactly the order statistics,
which include the classical median. By adding symmetry and/or idempotency to the previous results, we obtain the following corollaries.

\begin{corollary}
$A:\R^n\to\R$ is symmetric, nondecreasing (or continuous), and meaningful for the same input-output ordinal scales if and only if there exists
$k\in N$ such that $A={\rm OS}_k$.
\end{corollary}

\begin{corollary}
$A:\R^n\to\R$ is idempotent, nondecreasing (or continuous), and meaningful for the same input ordinal scales if and only if there exists
$\gamma\in\Gamma_N$ such that $A=L_{\gamma}$.
\end{corollary}

\begin{corollary}
$A:\R^n\to\R$ is symmetric, nondecreasing, and meaningful for the same input ordinal scales if and only if there exist $k\in N$ and a constant
or strictly increasing function $g:\R\to\R$ such that $A=g\circ {\rm OS}_k$.
\end{corollary}

\begin{corollary}
$A:\R^n\to\R$ is symmetric, continuous, and meaningful for the same input ordinal scales if and only if there exists $k\in N$ and a constant or
continuous and strictly monotonic function $g:\R\to\R$ such that $A=g\circ {\rm OS}_k$.
\end{corollary}

\section{Conclusion}

In this chapter we have discussed the most classical aggregation functions that are used in decision making. An appropriate classification of
these functions into a catalog can be better done through an axiomatic approach, which consists in listing a series of reasonable properties and
classifying or, better, characterizing the aggregation functions according to these properties.

With knowledge of the increasing need to define suitable aggregation functions fulfilling very precise conditions in various situations, it is not
surprising that such a catalog of aggregation functions, which is already huge, is constantly growing and remains an important topic of
research. We have only skimmed the surface of a still growing domain here.

\bibliographystyle{abbrv}   

\begin{thebibliography}{100}

\bibitem{AbrSte64}
M.~Abramowitz and I.~A. Stegun.
\newblock {\em Handbook of mathematical functions with formulas, graphs, and
  mathematical tables}, volume~55 of {\em National Bureau of Standards Applied
  Mathematics Series}.
\newblock For sale by the Superintendent of Documents, U.S. Government Printing
  Office, Washington, D.C., 1964.

\bibitem{Acz48}
J.~Acz{\'e}l.
\newblock On mean values.
\newblock {\em Bull. Amer. Math. Soc.}, 54:392--400, 1948.

\bibitem{Acz48b}
J.~Acz{\'e}l.
\newblock Sur les op\'erations d\'efinies pour nombres r\'eels.
\newblock {\em Bull. Soc. Math. France}, 76:59--64, 1948.

\bibitem{Acz66}
J.~Acz{\'e}l.
\newblock {\em Lectures on functional equations and their applications}.
\newblock Academic Press, New York, 1966.

\bibitem{AczAls87}
J.~Acz{\'e}l and C.~Alsina.
\newblock Synthesizing judgements: a functional equations approach.
\newblock {\em Math. Modelling}, 9(3-5):311--320, 1987.
\newblock The analytic hierarchy process.

\bibitem{AczDho89}
J.~Acz{\'e}l and J.~Dhombres.
\newblock {\em Functional equations in several variables}.
\newblock Cambridge University Press, Cambridge, 1989.
\newblock With applications to mathematics, information theory and to the
  natural and social sciences.

\bibitem{AczGroSch94}
J.~Acz{\'e}l, D.~Gronau, and J.~Schwaiger.
\newblock Increasing solutions of the homogeneity equation and of similar
  equations.
\newblock {\em J. Math. Anal. Appl.}, 182(2):436--464, 1994.

\bibitem{AczRob89}
J.~Acz{\'e}l and F.~S. Roberts.
\newblock On the possible merging functions.
\newblock {\em Math. Social Sci.}, 17(3):205--243, 1989.

\bibitem{AczRobRos86}
J.~Acz{\'e}l, F.~S. Roberts, and Z.~Rosenbaum.
\newblock On scientific laws without dimensional constants.
\newblock {\em J. Math. Anal. Appl.}, 119(1-2):389--416, 1986.

\bibitem{AlsMayTomTor93}
C.~Alsina, G.~Mayor, M.~S. Tom{\'a}s, and J.~Torrens.
\newblock A characterization of a class of aggregation functions.
\newblock {\em Fuzzy Sets and Systems}, 53(1):33--38, 1993.

\bibitem{AlsTriVal83}
C.~Alsina, E.~Trillas, and L.~Valverde.
\newblock On some logical connectives for fuzzy sets theory.
\newblock {\em J. Math. Anal. Appl.}, 93(1):15--26, 1983.

\bibitem{Ant98}
C.~Antoine.
\newblock {\em Les moyennes}, volume 3383 of {\em Que Sais-Je? [What Do I
  Know?]}.
\newblock Presses Universitaires de France, Paris, 1998.

\bibitem{BecBel65}
E.~F. Beckenbach and R.~Bellman.
\newblock {\em Inequalities}.
\newblock Second revised printing. Ergebnisse der Mathematik und ihrer
  Grenzgebiete. Neue Folge, Band 30. Springer-Verlag, New York, Inc., 1965.

\bibitem{Bem26}
G.~Bemporad.
\newblock Sul principio della media aritmetica. ({I}talian).
\newblock {\em Atti Accad. Naz. Lincei}, 3(6):87--91, 1926.

\bibitem{BerMor98}
L.~R. Berrone and J.~Moro.
\newblock Lagrangian means.
\newblock {\em Aequationes Math.}, 55(3):217--226, 1998.

\bibitem{BerMor00}
L.~R. Berrone and J.~Moro.
\newblock On means generated through the {C}auchy mean value theorem.
\newblock {\em Aequationes Math.}, 60(1-2):1--14, 2000.

\bibitem{Bir67}
G.~Birkhoff.
\newblock {\em Lattice theory}.
\newblock Third edition. American Mathematical Society Colloquium Publications,
  Vol. XXV. American Mathematical Society, Providence, R.I., 1967.

\bibitem{BulMitVas88}
P.~S. Bullen, D.~S. Mitrinovi{\'c}, and P.~M. Vasi{\'c}.
\newblock {\em Means and their inequalities}, volume~31 of {\em Mathematics and
  its Applications (East European Series)}.
\newblock D. Reidel Publishing Co., Dordrecht, 1988.
\newblock Translated and revised from the Serbo-Croatian.

\bibitem{Cau21}
A.~L. Cauchy.
\newblock {\em Cours d'analyse de l'Ecole Royale Polytechnique, Vol. I. Analyse
  alg\'ebrique}.
\newblock Debure, Paris, 1821.

\bibitem{Chi29}
O.~Chisini.
\newblock Sul concetto di media. ({I}talian).
\newblock {\em Periodico di matematiche}, 9(4):106--116, 1929.

\bibitem{Cho53}
G.~Choquet.
\newblock Theory of capacities.
\newblock {\em Ann. Inst. Fourier, Grenoble}, 5:131--295 (1955), 1953--1954.

\bibitem{Cli46}
A.~C. Climescu.
\newblock Sur l'\'equation fonctionnelle de l'associativit\'e.
\newblock {\em Bull. \'Ecole Polytech. Jassy [Bul. Politehn. Gh. Asachi. Ia\c
  si]}, 1:211--224, 1946.

\bibitem{CzoDre84}
E.~Czoga{\l}a and J.~Drewniak.
\newblock Associative monotonic operations in fuzzy set theory.
\newblock {\em Fuzzy Sets and Systems}, 12(3):249--269, 1984.

\bibitem{deCJor92}
L.~M. de~Campos and M.~Jorge.
\newblock Characterization and comparison of Sugeno and Choquet integrals.
\newblock {\em Fuzzy Sets and Systems}, 52(1):61--67, 1992.

\bibitem{deCLamMor91}
L.~M. de~Campos, M.~T. Lamata, and S.~Moral.
\newblock A unified approach to define fuzzy integrals.
\newblock {\em Fuzzy Sets and Systems}, 39(1):75--90, 1991.

\bibitem{deF31}
B.~de~{F}inetti.
\newblock Sul concetto di media. ({I}talian).
\newblock {\em Giorn. Ist. Ital. Attuari}, 2(3):369--396, 1931.

\bibitem{Del71}
C.~Dellacherie.
\newblock Quelques commentaires sur les prolongements de capacit\'es.
\newblock In {\em S\'eminaire de Probabilit\'es, V (Univ. Strasbourg, ann\'ee
  universitaire 1969-1970)}, pages 77--81. Lecture Notes in Math., Vol. 191.
  Springer, Berlin, 1971.

\bibitem{DowWer92}
J.~Dow and S.~R. d.~C. Werlang.
\newblock Uncertainty aversion, risk aversion, and the optimal choice of
  portfolio.
\newblock {\em Econometrica}, 60(1):197--204, 1992.

\bibitem{DowWer94}
J.~Dow and S.~R. d.~C. Werlang.
\newblock Nash equilibrium under {K}nightian uncertainty: breaking down
  backward induction.
\newblock {\em J. Econom. Theory}, 64(2):305--324, 1994.

\bibitem{DubMarPraRouSab01}
D.~Dubois, J.~L. Marichal, H.~Prade, M.~Roubens, and R.~Sabbadin.
\newblock The use of the discrete {S}ugeno integral in decision-making: a
  survey.
\newblock {\em Internat. J. Uncertain. Fuzziness Knowledge-Based Systems},
  9(5):539--561, 2001.

\bibitem{DubPra84}
D.~Dubois and H.~Prade.
\newblock Criteria aggregation and ranking of alternatives in the framework of
  fuzzy set theory.
\newblock In {\em Fuzzy sets and decision analysis}, volume~20 of {\em Stud.
  Management Sci.}, pages 209--240. North-Holland, Amsterdam, 1984.

\bibitem{DubPra85}
D.~Dubois and H.~Prade.
\newblock A review of fuzzy set aggregation connectives.
\newblock {\em Inform. Sci.}, 36(1-2):85--121, 1985.

\bibitem{DubPra86}
D.~Dubois and H.~Prade.
\newblock Weighted minimum and maximum operations in fuzzy set theory.
\newblock {\em Inform. Sci.}, 39(2):205--210, 1986.

\bibitem{DubPraTes88}
D.~Dubois, H.~Prade, and C.~Testemale.
\newblock Weighted fuzzy pattern matching.
\newblock {\em Fuzzy Sets and Systems}, 28(3):313--331, 1988.
\newblock Mathematical modelling.

\bibitem{Duj74}
J.~J. Dujmovi{\'c}.
\newblock Weighted conjunctive and disjunctive means and their application in
  system evaluation.
\newblock {\em Univ. Beograd. Publ. Elektrotehn. Fak. Ser. Mat. Fiz.},
  (461-497):147--158, 1974.

\bibitem{Duj75}
J.~J. Dujmovi{\'c}.
\newblock Extended continuous logic and the theory of complex criteria.
\newblock {\em Univ. Beograd. Publ. Elektrotehn. Fak. Ser. Mat. Fiz.},
  (498--541):197--216, 1975.

\bibitem{DycPed84}
H.~Dyckhoff and W.~Pedrycz.
\newblock Generalized means as model of compensative connectives.
\newblock {\em Fuzzy Sets and Systems}, 14(2):143--154, 1984.

\bibitem{Fau55}
W.~M. Faucett.
\newblock Compact semigroups irreducibly connected between two idempotents.
\newblock {\em Proc. Amer. Math. Soc.}, 6:741--747, 1955.

\bibitem{FisWak95}
P.~Fishburn and P.~Wakker.
\newblock The invention of the independence condition for preferences.
\newblock {\em Management Sciences}, 41(7):1130--1144, 1995.

\bibitem{Fod91}
J.~C. Fodor.
\newblock A remark on constructing {$t$}-norms.
\newblock {\em Fuzzy Sets and Systems}, 41(2):195--199, 1991.

\bibitem{Fod91b}
J.~C. Fodor.
\newblock Strict preference relations based on weak {$t$}-norms.
\newblock {\em Fuzzy Sets and Systems}, 43(3):327--336, 1991.
\newblock Aggregation and best choices of imprecise opinions (Brussels, 1989).

\bibitem{Fod96}
J.~C. Fodor.
\newblock An extension of {F}ung-{F}u's theorem.
\newblock {\em Internat. J. Uncertain. Fuzziness Knowledge-Based Systems},
  4(3):235--243, 1996.

\bibitem{FodRou94}
J.~C. Fodor and M.~Roubens.
\newblock {\em Fuzzy preference modelling and milticriteria decision support}.
\newblock Kluwer, Dordrecht, 1994.

\bibitem{Fro87}
V.~Frosini.
\newblock Averages.
\newblock In {\em Italian Contributions to the Methodology of Statistics},
  pages 1--17. Cleup, Padova, 1987.

\bibitem{FunFu75}
L.~W. Fung and K.~S. Fu.
\newblock An axiomatic approach to rational decision making in a fuzzy
  environment.
\newblock In {\em Fuzzy sets and their applications to cognitive and decision
  processes (Proc. U. S.-Japan Sem., Univ. Calif., Berkeley, Calif., 1974)},
  pages 227--256. Academic Press, New York, 1975.

\bibitem{Gil87}
I.~Gilboa.
\newblock Expected utility with purely subjective nonadditive probabilities.
\newblock {\em J. Math. Econom.}, 16(1):65--88, 1987.

\bibitem{Gra93}
M.~Grabisch.
\newblock On the use of fuzzy integral as a fuzzy connective.
\newblock In {\em Proc.\ of the Second IEEE International Conference on Fuzzy
  Systems}, pages 213--218, San Francisco, 1993.

\bibitem{Gra95}
M.~Grabisch.
\newblock Fuzzy integral in multicriteria decision making.
\newblock {\em Fuzzy Sets and Systems}, 69(3):279--298, 1995.

\bibitem{Gra95b}
M.~Grabisch.
\newblock On equivalence classes of fuzzy connectives: the case of fuzzy
  integrals.
\newblock {\em IEEE Trans. Fuzzy Syst.}, 3(1):96--109, 1995.

\bibitem{Gra96}
M.~Grabisch.
\newblock The application of fuzzy integrals in multicriteria decision making.
\newblock {\em European J. Oper. Res.}, 89(3):445--456, 1996.

\bibitem{Gra98}
M.~Grabisch.
\newblock Fuzzy integral as a flexible and interpretable tool of aggregation.
\newblock In {\em Aggregation and fusion of imperfect information}, volume~12
  of {\em Stud. Fuzziness Soft Comput.}, pages 51--72. Physica, Heidelberg,
  1998.

\bibitem{GraMarMesPap09}
M.~Grabisch, J.-L. Marichal, R.~Mesiar, and E.~Pap.
\newblock {\em Aggregation functions}.
\newblock Encyclopedia of Mathematics and its Applications. Cambridge
  University Press, Cambridge, UK, 2009.

\bibitem{GraMurSug00}
M.~Grabisch, T.~Murofushi, and M.~Sugeno, editors.
\newblock {\em Fuzzy measures and integrals, Theory and applications},
  volume~40 of {\em Studies in Fuzziness and Soft Computing}.
\newblock Physica-Verlag, Heidelberg, 2000.

\bibitem{GraNguWal95}
M.~Grabisch, H.~T. Nguyen, and E.~A. Walker.
\newblock {\em Fundamentals of uncertainty calculi with applications to fuzzy
  inference}, volume~30 of {\em Theory and Decision Library. Series B:
  Mathematical and Statistical Methods}.
\newblock Kluwer Academic Publishers, Dordrecht, 1995.

\bibitem{Gre87}
G.~H. Greco.
\newblock On {$L$}-fuzzy integrals of measurable functions.
\newblock {\em J. Math. Anal. Appl.}, 128(2):581--585, 1987.

\bibitem{GroJacSloTra98}
E.~Groes, J.~Jacobsen, B.~Sloth, and T.~Tran{\ae}s.
\newblock Axiomatic characterizations of the {C}hoquet integral.
\newblock {\em Econom. Theory}, 12(2):441--448, 1998.

\bibitem{Hoh82}
U.~H\"ohle.
\newblock Integration with respect to fuzzy measures.
\newblock In {\em Proc. IFAC Symposium on Theory and Applications of Digital
  Control}, pages 35--37, New Delhi, January 1982.

\bibitem{KanBya78}
A.~Kandel and W.~J. Byatt.
\newblock Fuzzy sets, fuzzy algebra, and fuzzy statistics.
\newblock {\em Proc. IEEE}, 66(12):1619--1639, December 1978.

\bibitem{KleMesPap00}
E.~P. Klement, R.~Mesiar, and E.~Pap.
\newblock {\em Triangular norms}, volume~8 of {\em Trends in Logic---Studia
  Logica Library}.
\newblock Kluwer Academic Publishers, Dordrecht, 2000.

\bibitem{Kol30}
A.~N. Kolmogoroff.
\newblock Sur la notion de la moyenne. ({F}rench).
\newblock {\em Atti Accad. Naz. Lincei}, 12(6):388--391, 1930.

\bibitem{Lin65}
C.-H. Ling.
\newblock Representation of associative functions.
\newblock {\em Publ. Math. Debrecen}, 12:189--212, 1965.

\bibitem{Luc59}
R.~D. Luce.
\newblock On the possible psychophysical laws.
\newblock {\em Psych. Rev.}, 66:81--95, 1959.

\bibitem{Mar98}
J.-L. Marichal.
\newblock {\em Aggregation operators for multicriteria decision aid}.
\newblock PhD thesis, Institute of Mathematics, University of Li\`ege, Li\`ege,
  Belgium, December 1998.

\bibitem{Mar00g}
J.-L. Marichal.
\newblock An axiomatic approach of the discrete choquet integral as a tool to
  aggregate interacting criteria.
\newblock {\em IEEE Trans. Fuzzy Syst.}, 8(6):800--807, 2000.

\bibitem{Mar00c}
J.-L. Marichal.
\newblock On {S}ugeno integral as an aggregation function.
\newblock {\em Fuzzy Sets and Systems}, 114(3):347--365, 2000.

\bibitem{Mar00b}
J.-L. Marichal.
\newblock On the associativity functional equation.
\newblock {\em Fuzzy Sets and Systems}, 114(3):381--389, 2000.

\bibitem{Mar01}
J.-L. Marichal.
\newblock An axiomatic approach of the discrete sugeno integral as a tool to
  aggregate interacting criteria in a qualitative framework.
\newblock {\em IEEE Trans. Fuzzy Syst.}, 9(1):164--172, 2001.

\bibitem{Mar02b}
J.-L. Marichal.
\newblock Aggregation of interacting criteria by means of the discrete choquet
  integral.
\newblock In {\em Aggregation operators: new trends and applications}, pages
  224--244. Physica, Heidelberg, 2002.

\bibitem{Mar02c}
J.-L. Marichal.
\newblock On order invariant synthesizing functions.
\newblock {\em J. Math. Psych.}, 46(6):661--676, 2002.

\bibitem{MarMatTou99}
J.-L. Marichal, P.~Mathonet, and E.~Tousset.
\newblock Characterization of some aggregation functions stable for positive
  linear transformations.
\newblock {\em Fuzzy Sets and Systems}, 102(2):293--314, 1999.

\bibitem{MarMesRuc}
J.-L. Marichal, R.~Mesiar, and T.~R\"uckschlossov\'a.
\newblock A complete description of comparison meaningful functions.
\newblock {\em Aequationes Math.}, 69(3):309--320, 2005.

\bibitem{MarRou93}
J.-L. Marichal and M.~Roubens.
\newblock Characterization of some stable aggregation functions.
\newblock In {\em Proc.\ 1st Int.\ Conf.\ on Industrial Engineering and
  Production Management (IEPM'93)}, pages 187--196, Mons, Belgium, June 1993.

\bibitem{MasMayTor99}
M.~Mas, G.~Mayor, and J.~Torrens.
\newblock t-operators.
\newblock {\em Internat. J. Uncertain. Fuzziness Knowledge-Based Systems},
  7(1):31--50, 1999.

\bibitem{MasMayTor02}
M.~Mas, G.~Mayor, and J.~Torrens.
\newblock The modularity condition for uninorms and t-operators.
\newblock {\em Fuzzy Sets and Systems}, 126(2):207--218, 2002.

\bibitem{MasMonTor01}
M.~Mas, M.~Monserrat, and J.~Torrens.
\newblock On left and right uninorms.
\newblock {\em Internat. J. Uncertain. Fuzziness Knowledge-Based Systems},
  9(4):491--507, 2001.

\bibitem{Mat99}
J.~Matkowski.
\newblock Invariant and complementary quasi-arithmetic means.
\newblock {\em Aequationes Math.}, 57(1):87--107, 1999.

\bibitem{MesRuc}
R.~Mesiar and T.~R\"uckschlossov\'a.
\newblock Characterization of invariant aggregation operators.
\newblock {\em Fuzzy Sets and Systems}, 142(1):63--73, 2004.

\bibitem{ModDubGraPra97}
F.~Modave, D.~Dubois, M.~Grabisch, and H.~Prade.
\newblock A choquet integral representation in multicriteria decision making.
\newblock In {\em AAAI Fall Symposium on Frontiers in Soft Computing and
  Decisions Systems}, pages 30--39, Boston, MA, November 7--9 1997.

\bibitem{ModGra98}
F.~Modave and M.~Grabisch.
\newblock Preference representation by the choquet integral: the
  commensurability hypothesis.
\newblock In {\em Proc.\ 7th Int.\ Conf.\ on Information Processing and
  Management of Uncertainty in Knowledge-Based Systems (IPMU'98)}, pages
  164--171, Paris, 1998.

\bibitem{MosShi57}
P.~S. Mostert and A.~L. Shields.
\newblock On the structure of semigroups on a compact manifold with boundary.
\newblock {\em Ann. of Math. (2)}, 65:117--143, 1957.

\bibitem{MurSug89}
T.~Murofushi and M.~Sugeno.
\newblock An interpretation of fuzzy measures and the {C}hoquet integral as an
  integral with respect to a fuzzy measure.
\newblock {\em Fuzzy Sets and Systems}, 29(2):201--227, 1989.

\bibitem{MurSug91}
T.~Murofushi and M.~Sugeno.
\newblock A theory of fuzzy measures: representations, the {C}hoquet integral,
  and null sets.
\newblock {\em J. Math. Anal. Appl.}, 159(2):532--549, 1991.

\bibitem{MurSug93}
T.~Murofushi and M.~Sugeno.
\newblock {Some quantities represented by the {C}hoquet integral}.
\newblock {\em Fuzzy Sets Syst.}, 56(2):229--235, 1993.

\bibitem{Nag30}
M.~Nagumo.
\newblock {\"U}ber eine klasse der mittelwerte. ({G}erman).
\newblock {\em Japanese Journ. of Math.}, 7:71--79, 1930.

\bibitem{NarMur02}
Y.~Narukawa and T.~Murofushi.
\newblock The $n$-step {C}hoquet integral on finite spaces.
\newblock In {\em 9th Int. Conf. on Information Processing and Management of
  Uncertainty in Knowledge-Based Systems (IPMU 2002)}, pages 539--543, Annecy,
  France, July 1-5 2002.

\bibitem{Orl81}
A.~I. Orlov.
\newblock The connection between mean values and the admissible transformations
  of scale.
\newblock {\em Math. Notes}, 30:774--778, 1981.

\bibitem{Ovc96}
S.~Ovchinnikov.
\newblock Means on ordered sets.
\newblock {\em Math. Social Sci.}, 32(1):39--56, 1996.

\bibitem{Ovc9798}
S.~Ovchinnikov.
\newblock Invariant functions on simple orders.
\newblock {\em Order}, 14(4):365--371, 1997/98.

\bibitem{RalRal97}
A.~L. Ralescu and D.~A. Ralescu.
\newblock Extensions of fuzzy aggregation.
\newblock {\em Fuzzy Sets and Systems}, 86(3):321--330, 1997.

\bibitem{RalSug96}
D.~A. Ralescu and M.~Sugeno.
\newblock Fuzzy integral representation.
\newblock {\em Fuzzy Sets and Systems}, 84(2):127--133, 1996.

\bibitem{Ric15}
U.~Ricci.
\newblock Confronti tra medie. ({I}talian).
\newblock {\em Giorn. Economisti e Rivista di Statistica}, 26:38--66, 1915.

\bibitem{Rob79}
F.~S. Roberts.
\newblock {\em Measurement theory}, volume~7 of {\em Encyclopedia of
  Mathematics and its Applications}.
\newblock Addison-Wesley Publishing Co., Reading, Mass., 1979.
\newblock With applications to decisionmaking, utility and the social sciences,
  Advanced Book Program.

\bibitem{Rob90}
F.~S. Roberts.
\newblock Merging relative scores.
\newblock {\em J. Math. Anal. Appl.}, 147(1):30--52, 1990.

\bibitem{Rob94}
F.~S. Roberts.
\newblock Limitations on conclusions using scales of measurement.
\newblock In {\em Operation Research and the Public Sector}, pages 621--671.
  Elsevier, Amsterdam, 1994.

\bibitem{San02}
W.~Sander.
\newblock Associative aggregation operators.
\newblock In {\em Aggregation operators}, volume~97 of {\em Stud. Fuzziness
  Soft Comput.}, pages 124--158. Physica, Heidelberg, 2002.

\bibitem{SarWak92}
R.~Sarin and P.~Wakker.
\newblock A simple axiomatization of nonadditive expected utility.
\newblock {\em Econometrica}, 60(6):1255--1272, 1992.

\bibitem{Sch86}
D.~Schmeidler.
\newblock Integral representation without additivity.
\newblock {\em Proc. Amer. Math. Soc.}, 97(2):255--261, 1986.

\bibitem{SchSkl61}
B.~Schweizer and A.~Sklar.
\newblock Associative functions and statistical triangle inequalities.
\newblock {\em Publ. Math. Debrecen}, 8:169--186, 1961.

\bibitem{SchSkl63}
B.~Schweizer and A.~Sklar.
\newblock Associative functions and abstract semigroups.
\newblock {\em Publ. Math. Debrecen}, 10:69--81, 1963.

\bibitem{SchSkl83}
B.~Schweizer and A.~Sklar.
\newblock {\em Probabilistic metric spaces}.
\newblock North-Holland Series in Probability and Applied Mathematics.
  North-Holland Publishing Co., New York, 1983.

\bibitem{Ste51}
S.~S. Stevens.
\newblock Mathematics, measurement, and psychophysics.
\newblock In {\em Handbook of Experimental Psychology}, pages 1--49. Wiley, New
  York, 1951.

\bibitem{Ste59}
S.~S. Stevens.
\newblock Measurement, psychophysics, and utility.
\newblock In {\em Measurement: Definitions and Theories}, pages 18--63. Wiley,
  New York, 1959.

\bibitem{Sug74}
M.~Sugeno.
\newblock {\em Theory of fuzzy integrals and its applications}.
\newblock PhD thesis, Tokyo Institute of Technology, Tokyo, 1974.

\bibitem{Sug77}
M.~Sugeno.
\newblock Fuzzy measures and fuzzy integrals---a survey.
\newblock In {\em Fuzzy automata and decision processes}, pages 89--102.
  North-Holland, New York, 1977.

\bibitem{Wak89}
P.~Wakker.
\newblock Continuous subjective expected utility with nonadditive
  probabilities.
\newblock {\em J. Math. Econom.}, 18(1):1--27, 1989.

\bibitem{Web83}
S.~Weber.
\newblock A general concept of fuzzy connectives, negations and implications
  based on {$t$}-norms and {$t$}-conorms.
\newblock {\em Fuzzy Sets and Systems}, 11(2):115--134, 1983.

\bibitem{YagKac97}
R.~Yager and J.~Kacprzyk, editors.
\newblock {\em The Ordered Weighted Averaging operators}.
\newblock Kluwer Academic Publishers, USA, 1997.
\newblock Theory and Applications.

\bibitem{Yag88}
R.~R. Yager.
\newblock On ordered weighted averaging aggregation operators in multicriteria
  decisionmaking.
\newblock {\em IEEE Trans. Systems Man Cybernet.}, 18(1):183--190, 1988.

\bibitem{YagRyb96}
R.~R. Yager and A.~Rybalov.
\newblock Uninorm aggregation operators.
\newblock {\em Fuzzy Sets and Systems}, 80(1):111--120, 1996.

\bibitem{Yan89}
E.~B. Yanovskaya.
\newblock Group choice rules in problems with interpersonal preference
  comparisons.
\newblock {\em Automat. Remote Control}, 50(6):822--830, 1989.

\bibitem{ZimZys80}
H.-J. Zimmermann and P.~Zysno.
\newblock Latent connectives in human decision making.
\newblock {\em Fuzzy Sets and Systems}, 4:37--51, 1980.

\end{thebibliography}

\end{document}